\documentclass{article}
\usepackage[utf8]{inputenc}

\pdfoutput=1
\usepackage{mathtools}
\usepackage{amsfonts}
\usepackage{amssymb}
\usepackage{centernot}
\usepackage{graphicx}
\usepackage{amsthm}
\usepackage{enumerate}

\usepackage{tikz}
\usepackage{tikz-cd}
\usetikzlibrary{decorations.markings}
\usetikzlibrary{arrows}
\usetikzlibrary{calc}
\providecommand{\U}[1]{\protect\rule{.1in}{.1in}}
\newtheorem{theorem}{Theorem}

\newtheorem{corollary}[theorem]{Corollary}

\newtheorem{definition}[theorem]{Definition}

\newtheorem{lemma}[theorem]{Lemma}

\newtheorem{proposition}[theorem]{Proposition}

\theoremstyle{remark}

\newcommand{\MG}{\textup{MG}}
\newcommand{\QMG}{\textup{QMG}}
\newcommand{\MQ}{\textup{MQ}}

\newcommand{\Aut}{\textup{Aut}}
\newcommand{\Dis}{\textup{Dis}}
\newcommand{\Fix}{\textup{Fix}}
\newcommand{\id}{\textup{id}}
\newcommand{\Qr}{Q_A^{\textup{red}}}
\newcommand{\Mr}{M_A^{\textup{red}}}
\newcommand{\fr}{\phi_L^{\textup{red}}}
\newcommand{\ann}{\textup{ann}}

\allowdisplaybreaks

\begin{document}

\title{Multivariate Alexander quandles, IV. The medial quandle of a link}
\author{Lorenzo Traldi\\Lafayette College\\Easton, PA 18042, USA\\traldil@lafayette.edu
}
\date{ }
\maketitle

\begin{abstract}
Joyce observed that the Alexander invariant and the medial quandle of a classical knot are equivalent to each other, as invariants. In the present paper, we discuss the rather complicated extension of Joyce's observation to several different medial quandles and reduced (one-variable) Alexander modules associated with classical links. The theme is that for links, medial quandles provide stronger invariants than reduced Alexander modules.


\end{abstract}

\section{Introduction}

If $L=K_1 \cup \dots \cup K_{\mu}$ is a classical, oriented link of $\mu$ components in $\mathbb S ^3$, then its (multivariate) Alexander module $M_A(L)$ is a module over the ring $\Lambda_{\mu}= \mathbb Z [t_1^{\pm 1}, \dots, t_{\mu}^{\pm 1}]$ of Laurent polynomials in $\mu$ variables, with integer coefficients. The Alexander module is often described by a presentation with generators and relations corresponding to arcs and crossings of a diagram of $L$; see Sec.\ \ref{defs} for details. 

A useful part of the Alexander module theory is the Crowell map, a module epimorphism $\phi_L:M_A(L) \to I_{\mu}$ introduced in \cite{C1}. Here $I_{\mu}$ is the augmentation ideal of $\Lambda_{\mu}$, i.e., the kernel of the augmentation map $\epsilon:\Lambda_{\mu} \to \mathbb Z$ given by $\epsilon(t_i^{\pm 1}) = 1 \thickspace \allowbreak \forall i \in \{1, \dots, \mu \}$. If $D$ is a diagram of $L$, then an element of $M_A(L)$ corresponding to an arc of $K_i$ in $D$ is mapped to $t_i-1$ by $\phi_L$.  We say that two links $L,L'$ are \emph{Crowell equivalent} if there is a module isomorphism $f:M_A(L) \to M_A(L')$ with $\phi_L=\phi_{L'} \circ f$. Following Rolfsen \cite{Ro}, we refer to $\ker \phi_L$ as the \emph{Alexander invariant} of $L$. The Alexander invariant corresponds to the first homology group of the universal abelian cover of $\mathbb S ^3 -L$, while the Alexander module corresponds to the first relative homology group of the covering space with respect to its fiber.

In the first paper of this series, $\phi_L$ was used to define an operation $\triangleright$ on $M_A(L)$. The operation $\triangleright$ defines quandle structures on various subsets of the Alexander module. One of these subsets yields $Q_A(L)$, the fundamental multivariate Alexander quandle of $L$. In the third paper of the series we completed a proof of the following.

\begin{theorem} (\cite{mvaq3, mvaq1})
\label{prequel1}
As an invariant of classical oriented links, the fundamental multivariate Alexander quandle $Q_A(L)$ (up to quandle isomorphism) is strictly stronger than $\phi_L$ (up to Crowell equivalence and permutation of component indices).
\end{theorem}

Here are three comments about Theorem \ref{prequel1}. 

(a) The theorem fails if index permutations are not allowed. For instance, let $L$ be a link whose Alexander polynomial is not symmetric with respect to a permutation of $\{1, \dots, \mu \}$, and let $L'$ be obtained from $L$ by applying that permutation to the component indices. Then $Q_A(L) \cong Q_A(L')$, as $Q_A(L)$ does not explicitly reflect component indices, but $M_A(L) \not \cong M_A(L')$. 

(b) The theorem implies that $Q_A(L)$ (up to quandle isomorphism) is also strictly stronger than both the Alexander invariant and the Alexander module (up to module isomorphism and permutation of $\{1, \dots, \mu \}$).

(c) The fundamental quandle $Q(L)$ is a strictly stronger link invariant than $Q_A(L)$.

The purpose of the present paper is to discuss the reduced (one-variable) version of the theory involved in Theorem \ref{prequel1}. Let $\Lambda=\mathbb Z [t^{\pm 1}]$ be the ring of Laurent polynomials in the variable $t$, with integer coefficients. If $\tau:\Lambda_{\mu} \to \Lambda$ is the homomorphism of rings with unity given by $\tau(t_i)=t\thickspace \allowbreak \forall i \in \{1, \dots, \mu \}$, then $\tau$ defines a $\Lambda_{\mu}$-module structure on $\Lambda$, with scalar multiplication given by $x \cdot y = \tau(x)y \thickspace \allowbreak \forall x \in \Lambda_{\mu}\thickspace \allowbreak \forall y \in \Lambda$. The \emph{reduced Alexander module} of $L$ is the tensor product
\[
M^{\textup{red}}_A(L)=M_A(L) \otimes _{\Lambda_{\mu}} \Lambda \textup{,}
\]
considered as a $\Lambda$-module via multiplication in the second factor. The tensor product of $\phi_L$ with the identity map of $\Lambda$ is a $\Lambda$-linear map 
\[
\phi_{\tau}:M^{\textup{red}}_A(L) \to I_{\mu} \otimes _{\Lambda_{\mu}} \Lambda.
\]

\begin{definition}
\label{taueq}
Two links $L,L'$ are \emph{$\phi_{\tau}$-equivalent} if there is a $\Lambda$-module isomorphism $f:M^{\textup{red}}_A(L) \to M^{\textup{red}}_A(L')$ that is compatible with the $\phi_{\tau}$ maps of $L$ and $L'$, i.e., $\phi_{\tau}=\phi'_{\tau} \circ f$.
\end{definition}

\begin{definition}
\label{qred}
If $L$ is a classical link, let
\[
Q_A^{\textup{red}}(L)=\{x \otimes 1 \mid x \in Q_A(L)\} \subset M_A^{\textup{red}}(L).
\]
\end{definition}

It is not hard to verify that the quandle operation $\triangleright$ of $Q_A(L)$ defines a quandle structure on $Q_A^{\textup{red}}(L)$ in a natural way: $(x \otimes 1) \triangleright (y \otimes 1)=(x \triangleright y) \otimes 1$. Also, $Q_A^{\textup{red}}(L)$ is a subquandle of the standard Alexander quandle on the $\Lambda$-module $M_A^{\textup{red}}(L)$. In fact, $Q_A^{\textup{red}}(L)$ is an invariant subquandle of $M_A^{\textup{red}}(L)$, in this sense: if $L$ and $L'$ are ambient isotopic oriented links, then there is an isomorphism $M_A^{\textup{red}}(L) \cong M_A^{\textup{red}}(L')$ that maps $Q_A^{\textup{red}}(L)$ isomorphically onto $Q_A^{\textup{red}}(L')$.

In Section \ref{proof1} we verify the following.

\begin{theorem}
\label{main1}
As an invariant of classical links, $Q_A^{\textup{red}}(L)$ (up to quandle isomorphism) is equivalent to $\phi_{\tau}$ (up to $\phi_{\tau}$-equivalence and permutation of component indices).
\end{theorem}

Modified versions of the earlier comments (a), (b), (c) hold for Theorem \ref{main1}. The easiest one to state is (a): like Theorem \ref{prequel1}, Theorem \ref{main1} fails if index permutations are disallowed. See Sec.\ \ref{hopfex} for details.

In order to state the reduced version of comment (b), it is convenient to let $\eta:I_{\mu} \otimes _{\Lambda_{\mu}} \Lambda \to \Lambda$ be the map with $\eta((t_i-1) \otimes 1) = 1 \thickspace \allowbreak \forall i \in \{1, \dots, \mu \}$. The composition $\fr = \eta \circ \phi_\tau:\Mr(L) \to \Lambda$ appears in the reduced version of Crowell's link module sequence, i.e.\ the homology sequence of the total linking number cover of $\mathbb S ^3 -L$. We call $\ker \fr$ the \emph{reduced Alexander invariant} of $L$. For more information regarding the properties of reduced link module sequences, we refer to Hillman \cite[Sec.\ 5.4]{H}.

We now have three $\Lambda$-modules associated with $L$: the reduced Alexander module $\Mr(L)$, the reduced Alexander invariant $\ker \fr$, and $\ker \phi_\tau$, for which we do not have a special name. The situation  may seem complicated, but it turns out that one module determines the other two.

\begin{proposition}
The $\Lambda$-module $\ker \fr$ determines both $\Mr(L)$ and $\ker \phi_\tau$, up to isomorphism: $\Mr(L) \cong \ker \fr \oplus \Lambda $ and $ \ker \phi_\tau = (t-1) \cdot \ker \fr$.
\end{proposition}
\begin{proof}
The epimorphism $\fr:\Mr \to \Lambda$ must split, so $\Mr(L) \cong \ker \fr \oplus \Lambda$. The equality $\ker \phi_\tau = (t-1) \cdot \ker \fr$ is not so obvious; see Proposition \ref{kerphi}.
\end{proof}

Here is the reduced version of comment (b).

\begin{theorem}
\label{maincor}
In general, the quandle $Q_A^{\textup{red}}(L)$ (up to quandle isomorphism) is a strictly stronger link invariant than the module $\ker \fr$. For knots, though, $Q_A^{\textup{red}}(L)$ and $\ker \fr$ are equivalent invariants.
\end{theorem}

In addition to being the image of $Q_A(L)$ in $M_A^{\textup{red}}(L)$, $Q_A^{\textup{red}}(L)$ is also the image in $M_A^{\textup{red}}(L)$ of a link invariant introduced by Joyce \cite{J}, namely, the fundamental medial quandle $\textup{MQ}(L)$. (Joyce denoted this quandle $\textup{AbQ}(L)$ rather than $\textup{MQ}(L)$, and he called it the ``abelian link quandle'' of $L$.) Joyce proved that when $\mu=1$, $\textup{MQ}(L)$ and $\ker \phi_L$ are equivalent invariants \cite[Sec.\ 17]{J}. This property is reflected in the reduced version of comment (c):
 
\begin{theorem}
\label{main2} For knots, $\MQ(L)$ and $\Qr(L)$ are isomorphic quandles. In general, though, $\textup{MQ}(L)$ is a strictly stronger link invariant than $\Qr(L)$. 
\end{theorem}

We should mention that early versions of the present paper, posted on the arxiv, included the incorrect assertion that $\Qr(L)$ and $\MQ(L)$ are always isomorphic. We are grateful to Kyle Miller for helping us understand the mistake.

Here is an outline of our discussion.  In Sec.\ \ref{defs}, we present some basic properties of Alexander modules and Crowell maps. In Sec.\ \ref{proof1}, we discuss several quandles associated with Alexander modules, and prove Theorem \ref{main1}. Theorem \ref{maincor} is proven in Sec.\ \ref{hopfex}. 

In the rest of the paper, we do more than just prove Theorem \ref{main2}; we try to provide as much insight as we can into the quandles $\Qr(L)$ and $\MQ(L)$, and their connections with the $\Lambda$-module $\Mr(L)$. Sec.\ \ref{struc} is a brief account of some of the basic theory of general quandles and medial quandles; most of the material is drawn from the work of Jedli\v{c}ka, Pilitowska, Stanovsk\'{y} and Zamojska-Dzienio \cite{JPSZ1, JPSZ2}. In Sec.\ \ref{linkq}, we discuss Joyce's description of $\MQ(L)$ as a quandle ``augmented'' by a group \cite{J}; we denote this group $\MG(L)$. We also verify two assertions of Theorem \ref{main2}: $\MQ(L)$ determines $\Qr(L)$, and if $\mu=1$, then $\MQ(L) \cong \Qr(L)$. In Sec.\ \ref{twoproof}, we complete the proof of Theorem \ref{main2} by providing examples distinguished by $\MQ(L)$ but not by $\Qr(L)$. In the last two sections of the paper, we show that $Q_A^{\textup{red}}(L)$ is isomorphic to a quandle contained in the group $\MG(L)$.

\section{Alexander Modules and Crowell Maps}
\label{defs}

We follow the usual conventions for diagrams of classical links. A diagram $D$ consists of piecewise smooth closed curves in the plane, whose only (self-) intersections are \emph{crossings}, i.e., transverse double points. The set of crossings in $D$ is denoted $C(D)$. At each crossing, two short segments are removed, to indicate which of the intersecting curves is the underpasser. Removing these short segments cuts the curves into separate parts, the \emph{arcs} of $D$. The set of arcs is denoted $A(D)$. If $D$ is a diagram of $L=K_1 \cup \dots \cup K_{\mu}$, then there is a function $\kappa_D:A(D) \to \{1, \dots, \mu\}$, with $\kappa_D(a)=i$ if $a$ is part of the image of $K_i$ in $D$.
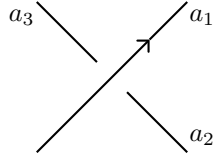
\begin{figure} [bth]
\centering
\begin{tikzpicture} [>=angle 90]
\draw [thick] (1,1) -- (0.5,0.5);
\draw [thick] [<-] (0.5,0.5) -- (-1,-1);
\draw [thick] (-1,1) -- (-.2,0.2);
\draw [thick] (0.2,-0.2) -- (1,-1);
\node at (1.2,0.8) {$a_1$};
\node at (-1.2,0.8) {$a_3$};
\node at (1.2,-0.8) {$a_2$};
\end{tikzpicture}
\caption{A crossing.}
\label{crossfig}
\end{figure}

Let $D$ be a link diagram, and let $\Lambda_{\mu}^{A(D)}$ and $\Lambda_{\mu}^{C(D)}$ be the free $\Lambda_{\mu}$-modules on the sets $A(D)$ and $C(D)$. There is a $\Lambda_{\mu}$-linear map $\rho_D:\Lambda_{\mu}^{C(D)} \to \Lambda_{\mu}^{A(D)}$ given by
\[
\rho_D(c)=(1-t_{\kappa_D(a_2)})a_1+t_{\kappa_D(a_1)}a_2-a_3
\]
whenever $c \in C(D)$ is a crossing of $D$ as indicated in Fig.\ \ref{crossfig}, and there is an exact sequence 
\begin{equation*}
\Lambda_{\mu}^{C(D)} \xrightarrow{\rho_D} \Lambda_{\mu}^{A(D)} \xrightarrow{\gamma_D} M_A(L) \to 0 \textup{.}
\end{equation*}
The Crowell map $\phi_L:M_A(L) \to I_{\mu}$ has $\phi_L\gamma_D(a)=t_{\kappa_D(a)}-1 \thickspace \allowbreak \forall a \in A(D)$. 

The reduced Alexander module of $L$ is $M^{\textup{red}}_A(L)=M_A(L) \otimes_{\Lambda{\mu}} \Lambda$. It is equivalent to say that $M_A^{\textup{red}}(L)$ is the quotient $M_A(L)/(J \cdot M_A(L))$, where $J$ is the ideal of $\Lambda_{\mu}$ generated by the elements $t_i-t_j$. (An isomorphism between $M_A(L)/(J \cdot M_A(L))$ and $M_A(L) \otimes _{\Lambda_{\mu}} \Lambda$ is given by $x+ J \cdot M_A(L) \leftrightarrow x \otimes 1$.) Alternatively, $M_A^{\textup{red}}(L)$ can be obtained by setting all $t_i$ equal to $t$ in any description of $M_A(L)$ (e.g., the description in the preceding paragraph).

Every element of $M^{\textup{red}}_A(L)$ is of the form $x \otimes 1$ for some $x \in M_A(L)$, for the following reasons: (a) if $x \in M_A(L)$ and $\lambda \in \Lambda_{\mu}$ then $x \otimes \tau(\lambda)=(\lambda 
\cdot x) \otimes 1$, and (b) if $x_1,x_2 \in M_A(L)$ then $x_1 \otimes 1 + x_2 \otimes 1 = (x_1 + x_2) \otimes 1$. It follows that formulas involving elements of $M^{\textup{red}}_A(L)$ can be specified using elements of the form $x \otimes 1$, with $x \in M_A(L)$. For instance, the map $\phi_{\tau}:M^{\textup{red}}_A(L) \to I_{\mu} \otimes_{\Lambda{\mu}} \Lambda$ is defined by $\phi_{\tau}(x \otimes 1)=\phi_L(x) \otimes 1 \thickspace \allowbreak \forall x \in M_A(L)$, and if $D$ is a diagram of $L$, then $\phi_{\tau}(\gamma_D(a) \otimes 1)=(t_{\kappa_D(a)}-1) \otimes 1 \thickspace \allowbreak \forall a \in A(D)$. 

When thinking about the map $\phi_{\tau}$, it is helpful to have in mind an explicit description of the $\Lambda$-module $I_{\mu} \otimes_{\Lambda{\mu}} \Lambda$. Let $\epsilon:\Lambda \to \mathbb Z$ be the augmentation map given by $\epsilon(t)=1$, and let $\mathbb Z_{\epsilon}$ be the $\Lambda$-module obtained from $\mathbb Z$ using the scalar multiplication given by $\lambda \cdot n = \epsilon(\lambda)n \thickspace \allowbreak \forall \lambda \in \Lambda \thickspace \allowbreak \forall n \in \mathbb Z$. Notice that $\mathbb Z _ \epsilon \cong \Lambda/(t-1)$, where $(t-1) = \ker \epsilon$ is the augmentation ideal of $\Lambda$.
\begin{lemma}
\label{tenstruc}
As a $\Lambda$-module, $I_{\mu} \otimes_{\Lambda{\mu}} \Lambda$ is isomorphic to
\[
\Lambda \oplus \underbrace{ \mathbb Z _{\epsilon} \oplus \cdots \oplus \mathbb Z _ \epsilon}_{\mu-1} \textup{,}
\]
with the $\Lambda$ summand generated by $(t_1-1) \otimes 1$ and the $\mathbb Z _ \epsilon$ summands generated by $(t_2-t_1) \otimes 1, \dots, (t_{\mu}-t_1) \otimes 1 $.
\end{lemma}
\begin{proof}
It is well known that as a $\Lambda_{\mu}$-module, $I_{\mu}$ is generated by the elements $t_1-1, \dots, t_{\mu}-1$, subject to the defining relations $(t_i-1) \cdot (t_j-1)=(t_j-1) \cdot (t_i-1)\thickspace \allowbreak \forall i,j$. (See \cite[p.\ 71]{H}, for instance.) It follows that as a $\Lambda$-module, $I_{\mu} \otimes_{\Lambda{\mu}} \Lambda$ is generated by the elements $(t_1-1) \otimes 1, \dots, (t_{\mu}-1) \otimes 1$, subject to the defining relations 
\[
(t-1) \cdot ((t_j-1) \otimes 1)=(t_j-1) \otimes (t-1)=(t_j-1) \otimes \tau(t_i-1)=((t_j-1)(t_i-1))\otimes 1
\]
\[
=((t_i-1)(t_j-1))\otimes 1 = (t_i-1)\otimes \tau(t_j-1)=(t_i-1) \otimes (t-1)=(t-1) \cdot ((t_i-1) \otimes 1)
\]
for all values of $i$ and $j$. Equivalently, $I_{\mu} \otimes_{\Lambda{\mu}} \Lambda$ is generated by $(t_1-1) \otimes 1$ and the $\mu -1$ elements
\[
(t_2-1) \otimes 1 - (t_1-1) \otimes 1, \dots, (t_{\mu}-1) \otimes 1 - (t_1-1) \otimes 1 \textup{,}
\]
subject to the defining relations $(t-1) \cdot ((t_i-1) \otimes 1 - (t_1-1) \otimes 1) = 0 \thickspace \allowbreak \forall i \in \{2, \dots, \mu \}$. As $(t_i-1) \otimes 1 - (t_1-1) \otimes 1=(t_i-t_1) \otimes 1\thickspace \allowbreak \forall i \in \{2, \dots, \mu \}$, the result follows. \end{proof}

The map $\eta:I_{\mu} \otimes _{\Lambda_{\mu}} \Lambda \to \Lambda$ was defined in the introduction by the formula $\eta((t_i-1) \otimes 1) = 1 \thickspace \allowbreak \forall i \in \{1, \dots, \mu \}$. Taking the isomorphism of Lemma 8 into account, $\eta:\Lambda \oplus (\mathbb Z _\epsilon )^{\mu -1} \to \Lambda$ is simply the projection onto the first coordinate of the direct sum. Therefore $\fr = \eta \circ \phi_\tau:\Mr(L) \to \Lambda$ is the first coordinate of $\phi_\tau:\Mr(L) \to \Lambda \oplus (\mathbb Z _\epsilon )^{\mu -1}$. As $\fr(\gamma_D(a) \otimes 1)=1\thickspace \allowbreak \forall a \in A(D)$, it is easy to see that $\ker \fr$ is the submodule of $\Mr(L)$ generated by the elements $(\gamma_D(a)-\gamma_D(a')) \otimes 1$ with $a,a' \in A(D)$. It is not much harder to describe $\ker \phi_\tau$.

\begin{proposition}
\label{kerphi}
The kernel of $\phi_{\tau}$ is $(t-1) \cdot \ker \fr$.
\end{proposition}
\begin{proof} Identify $I_{\mu} \otimes_{\Lambda{\mu}} \Lambda$ with $\Lambda \oplus (\mathbb Z _ \epsilon)^{\mu-1}$ using the isomorphism of Lemma \ref{tenstruc}. If $x \in \ker \fr$, then $\phi_\tau(x) = (0,n_1, \dots n_{\mu-1}) \in \Lambda \oplus (\mathbb Z _ \epsilon)^{\mu-1}$ for some integers $n_1, \dots n_{\mu-1}$. As $t \cdot n_i = n_i$ in $\mathbb Z _ \epsilon$ for each $i$, it follows that $\phi_\tau((t-1)x)=(t-1)\phi_\tau(x)=0$. Thus $(t-1) \cdot \ker \fr \subseteq \ker \phi_\tau$.

To verify the opposite inclusion, suppose $a_1,\dots,a_n \in A(D)$, $\lambda_1,\dots,\lambda_n \in \Lambda$, and 
\begin{equation}
\label{xsum}
x = \sum _{j=1}^n \lambda_j \cdot  (\gamma_D(a_j)\otimes 1) \in \ker \phi_{\tau}.
\end{equation}
We claim that $x \in (t-1) \cdot \ker \fr$.

Every coordinate of $\phi_{\tau}(x)$ in $\Lambda \oplus (\mathbb Z _ \epsilon)^{\mu-1}$ is $0$, so
\begin{align*}
& \textup{(i)  } \sum_{j=1}^n \lambda_j=0 \textup{  and}\\
& \textup{(ii)  for each } i \in \{2,\dots,\mu\} \textup{, } 
\sum _{\kappa_D(a_j)=i} \epsilon (\lambda_j )=0 .
\end{align*}
Notice that property (i) implies that property (ii) holds also when $i=1$. 

Suppose first that $\epsilon(\lambda_1),\dots,\epsilon(\lambda_n)$ are all $0$. The kernel of $\epsilon$ is the principal ideal $(t-1)$ of $\Lambda$, so for each $j$, there is a $\lambda'_j \in \Lambda$ with $\lambda_j=\lambda'_j \cdot (t-1)$. Property (i) implies that $\sum_{j=1}^n \lambda'_j=0$, so if $a^*$ is any fixed element of $A(D)$,
\[
x = (t-1)\sum_{j=1}^n \lambda'_j \cdot  (\gamma_D(a_j)\otimes 1) - (t-1) \bigg(\sum_{j=1}^n \lambda'_j \bigg) \cdot  (\gamma_D(a^*)\otimes 1)
\]
\[
= \sum_{j=1}^n \lambda'_j \cdot (t-1) \cdot  ((\gamma_D(a_j) - \gamma_D(a^*)) \otimes 1).
\]
This satisfies the claim.

Now, suppose at least one of $\epsilon(\lambda_1),\dots,\epsilon(\lambda_n)$ is not $0$. For convenience, introduce a new summand $0 \cdot (\gamma_D(a) \otimes 1)$ for each $ a \in A(D)$, and collect all the appearances of each $a_j$ into one summand, so that every $a \in A(D)$ appears precisely once in (\ref{xsum}). If all values of $\epsilon(\lambda_j)$ are now $0$, the earlier argument applies. Otherwise, re-index the elements of $A(D)$ so that $\epsilon(\lambda_1) \neq 0$. Let $\kappa=\kappa_D(a_1)$. Re-index the elements of $A(D)$ so that for some $k \in \{1, \dots, n \}$, $a_1, \dots, a_k$ are the arcs of $D$ with $\kappa_D(a_i)=\kappa$, and the arcs $a_1, \dots, a_k,a_1$ are encountered in this order as we walk along the image of $K_{\kappa}$ in $D$. Notice that according to property (ii) above, $\epsilon(\lambda_1) \neq 0$ implies that $k>1$, so $a_1, \dots, a_k$ are all distinct.

For $1 \leq i < k$, let $a'_i$ be the overpassing arc at the crossing $c_i$ of $D$ that separates $a_i$ from $a_{i+1}$. Then depending on the orientation of $a'_i$, one of these two formulas is equal to $0$ in $M_A^{\textup{red}}(L)$.
\[
\gamma_D\rho_D(c_i) \otimes 1 = (1-t)(a'_i \otimes 1)+t(a_i\otimes 1)-a_{i+1}\otimes 1
\]
\[
-\gamma_D\rho_D(c_i) \otimes 1 = -(1-t)(a'_i \otimes 1)-t(a_{i+1}\otimes 1)+a_i\otimes 1
\]
Let $0_i$ denote one of the two displayed formulas that does equal $0$ in $M_A^{\textup{red}}(L)$. Notice that if we add $0_i$ to the sum (\ref{xsum}), the only effect on the values of $\epsilon(\lambda_1),\dots,\epsilon(\lambda_n)$ is to add $1$ to the value of $\epsilon(\lambda_i)$, and add $-1$ to the value of $\epsilon(\lambda_{i+1})$. Of course if we add $-0_i$ instead, we produce the opposite effects.

It follows that by repeatedly adding $\pm 0_1$ to the sum in (\ref{xsum}), we can obtain a sum still equal to $x$, in which $\epsilon(\lambda_1)$ is $0$. Doing the same thing for $i=2, \dots, k-1$, we obtain a sum still equal to $x$, in which $\epsilon(\lambda_1), \dots, \epsilon(\lambda_{k-1})$ are all $0$. Property (ii) then implies that $\epsilon(\lambda_k)$ is $0$ too, so every arc $a_j \in A(D)$ with $\kappa_D(a_j)=\kappa$ has $\epsilon(\lambda_j)=0$.

Repeating this argument for each component of $L$ that has some arc $a_j$ with $\epsilon(\lambda_j) \neq 0$, we ultimately obtain a sum (\ref{xsum}) equal to $x$ in which $\epsilon(\lambda_1),\dots,\epsilon(\lambda_n)$ are all $0$. Then the earlier argument tells us that the claim holds for $x$. \end{proof}

We end this section with a well-known property of the Alexander invariants of knots.

\begin{corollary}
\label{knotker}
Suppose $\mu=1$. Then scalar multiplication by $t-1$ defines an automorphism of $\ker \fr$ as a $\Lambda$-module.
\end{corollary}
\begin{proof}
Of course, scalar multiplication by any element of $\Lambda$ defines an endomorphism of any $\Lambda$-module. As $\mu=1$, the map $\eta:I_{\mu} \otimes _{\Lambda_{\mu}} \Lambda \to \Lambda$ is an isomorphism, so $\ker \phi_\tau= \ker (\eta \phi_\tau) = \ker \fr$; hence Proposition \ref{kerphi} tells us that $(t-1) \cdot \ker \fr = \ker \fr$. Therefore, scalar multiplication by $t-1$ defines a surjective endomorphism of $\ker \fr$. As a submodule of the finitely generated module $\Mr(L)$ over the Noetherian ring $\Lambda$, $\ker \fr$ is Noetherian; therefore a surjective endomorphism of $\ker \fr$ must be an automorphism. \end{proof}

\section{Theorem \ref{main1}}
\label{proof1}

Recall the definition of a quandle:

\begin{definition}
\label{qdef}
A \emph{quandle} is a set $Q$ equipped with a binary operation $\triangleright$, which satisfies the following properties.
\begin{enumerate}
\item $x\triangleright x=x \thickspace \allowbreak \forall x \in Q$.
\item For each $y \in Q$, the formula $\beta_y(x)=x \triangleright y$ defines a permutation $\beta_y$ of $Q$.
\item $(x\triangleright y) \triangleright z=(x\triangleright z) \triangleright (y\triangleright z) \thickspace \allowbreak \forall x,y,z \in Q$.
\end{enumerate}
\end{definition}
There is a traditional way to associate a quandle to any $\Lambda$-module, which was mentioned by both Joyce and Matveev when they introduced quandles as link invariants \cite{J, M}.

\begin{proposition} (\cite{J, M})
\label{standardq}
If $M$ is a $\Lambda$-module, then the operations $x \triangleright y=tx+(1-t)y$ and $x \triangleright^{-1} y=t^{-1} x+(1-t^{-1})y$ define a quandle structure on $M$.
\end{proposition}

The quandles described in Proposition \ref{standardq} are called \emph{Alexander quandles} in the literature. In order to distinguish them from other quandles associated with Alexander modules, we refer to them as \emph{standard} Alexander quandles. Notice that every standard Alexander quandle is a whole $\Lambda$-module. Also, if $M$ is a standard Alexander quandle and $w,x,y,z \in M$, then
\[
(w \triangleright x)\triangleright (y \triangleright z) = t^2w + t(1-t)x +(1-t)ty+(1-t)^2z
\]
\[
= t^2w + t(1-t)y +(1-t)tx+(1-t)^2z=(w \triangleright y) \triangleright (x \triangleright z).
\]
That is, all standard Alexander quandles satisfy the \emph{medial} property $(w \triangleright x)\triangleright (y \triangleright z)=(w \triangleright y)\triangleright (x \triangleright z)$. It follows that all subquandles of standard Alexander quandles are medial, too.

In \cite{mvaq1} we introduced the operation $x \triangleright y = (\phi_L(y)+1)x - \phi_L(x) y$ on $M_A(L)$, and showed that $\triangleright$ defines a quandle structure on the subset 
\[
U(L) = \{x \in M_A(L) \mid \phi_L(x)+1 \textup{ is a unit of } \Lambda_{\mu} \}.
\]
The subquandle of $U(L)$ generated by $\gamma_D(A(D))$ is the fundamental multivariate Alexander quandle, $Q_A(L)$. 
The following result was proven in \cite{mvaq1}.

\begin{theorem} (\cite{mvaq1})
\label{mqinv}
If $L$ and $L'$ are equivalent links with diagrams $D$ and $D'$, then there is an isomorphism $f:M_A(L) \to M_A(L')$ which maps the quandle $U(L)$ isomorphically onto $U(L')$, maps the quandle $Q_A(L)$ isomorphically onto $Q_A(L')$, and is compatible with the Crowell maps of $L$ and $L'$, i.e., $\phi_L=\phi_{L'} \circ f$.
\end{theorem}

Multivariate Alexander quandles differ from standard Alexander quandles in several regards. For one thing, $Q_A(L)$ corresponds to a proper subset of the $\Lambda_{\mu}$-module $M_A(L)$, not a whole $\Lambda$-module. For another, $Q_A(L)$ is not a medial quandle, in general. 

In the introduction, we defined $Q_A^{\textup{red}}(L)$ to be $\{ x \otimes 1 \mid x \in Q_A(L)\} \subset M_A^{\textup{red}}(L)$, and stated that it is a quandle under the operation $(x \otimes 1) \triangleright (y \otimes 1)=(x \triangleright y) \otimes 1$. Here is an equivalent description.

\begin{proposition}
\label{standardsub}
Let $D$ be a diagram of a link $L$. Then $Q_A^{\textup{red}}(L)$ is the subquandle of the standard Alexander quandle on $M^{\textup{red}}_A(L)$ generated by the elements $\gamma_D(a) \otimes 1$, $a \in A(D)$.
\end{proposition}
\begin{proof}
If $x,y \in Q_A(L)$ then $\phi_L(x),\phi_L(y) \in \{t_1-1, \dots, t_{\mu}-1\}$, so
\[
(x \otimes 1) \triangleright (y \otimes 1) = (x \triangleright y) \otimes 1 = ((\phi_L(y)+1)x - \phi_L(x) y) \otimes 1 
\]
\[
= ((\phi_L(y)+1)x)\otimes 1 - (\phi_L(x) y) \otimes 1
= x \otimes \tau(\phi_L(y)+1)-y \otimes \tau(\phi_L(y))
\]
\[
=x \otimes t - y \otimes (t-1) = t \cdot (x \otimes 1) + (1-t) \cdot (y \otimes 1) \textup{.}
\]
This equals $(x \otimes 1) \triangleright (y \otimes 1)$ in the standard Alexander quandle on $M^{\textup{red}}_A(L)$. 

The result follows, as $Q_A(L)$ is generated by $\gamma_D(A(D))$.  \end{proof}

Having both descriptions of $Q_A^{\textup{red}}(L)$ is convenient because it makes it unnecessary to provide new proofs for many properties of $Q_A^{\textup{red}}(L)$. Instead, we can simply refer to established properties of $Q_A(L)$ and $M^{\textup{red}}_A(L)$. For instance, Propositions \ref{standardq} and \ref{standardsub} tell us that $Q_A^{\textup{red}}(L)$ is indeed a quandle. 

Here are some other properties of $Q_A^{\textup{red}}(L)$.

\begin{theorem} 
Suppose $L$ and $L'$ are oriented links of the same link type, with diagrams $D$ and $D'$, and associated maps
\[
\phi_{\tau}:M_A^{\textup{red}}(L) \to I_{\mu} \otimes_{\Lambda_{\mu}} \Lambda \quad \textup{  and  } \quad \phi'_{\tau}:M_A^{\textup{red}}(L') \to I_{\mu} \otimes_{\Lambda_{\mu}} \Lambda.
\]
Then there is an isomorphism $M_A^{\textup{red}}(L) \cong M_A^{\textup{red}}(L')$, which maps $Q_A^{\textup{red}}(L)$ isomorphically onto $Q_A^{\textup{red}}(L')$, and has $\phi_{\tau}=\phi'_{\tau} \circ f$.
\end{theorem}
\begin{proof}
The result follows from Theorem \ref{mqinv}, using the right exactness of tensor products. A direct proof can be obtained by setting all $t_i$ equal to $t$ in the discussion of \cite[Sec.\ 3]{mvaq1}.
\end{proof}

Recall that an \emph{orbit} in a quandle $Q$ is an equivalence class under the equivalence relation generated by $x \sim x \triangleright^{-1} y \sim x \triangleright y \thickspace \allowbreak \forall x,y \in Q$.

\begin{theorem}
\label{qpres}
Let $L=K_1 \cup \dots \cup K_{\mu}$ be a link with a diagram $D$, and let $Q$ be $Q_A^{\textup{red}}(L)$, or $Q_A(L)$, or the fundamental quandle $Q(L)$. Then there are surjective functions $\kappa_D:Q \to \{1, \dots, \mu\}$ and $\widehat{\sigma}_{\tau}:\Lambda^Q \to M_A^{\textup{red}}(L)$ with the following properties.
\begin{enumerate}
    \item The orbits of $Q$ are the sets $\kappa_D^{-1}(\{1\}),\dots,\kappa_D^{-1}(\{\mu\})$. Moreover, after an appropriate permutation of the indices of $K_1, \dots, K_\mu$, it will be true that $\kappa_D(q)=\kappa_D(a)$ whenever $q \in Q$ corresponds to an arc $a \in A(D)$.
    \item The restriction $\sigma_{\tau}=\widehat{\sigma}_{\tau}|Q$ is a quandle homomorphism onto $Q_A^{\textup{red}}(L)$, and if $a \in A(D)$ then the image under $\sigma_{\tau}$ of the element of $Q$ corresponding to $a$ is $\gamma_D(a) \otimes 1$. In particular, if $Q=Q^{\textup{red}}_A(L)$ then $\sigma_{\tau}$ is the identity map of $Q$.
    \item The map $\widehat{\sigma}_{\tau}$ is $\Lambda$-linear, and $\ker \widehat{\sigma}_{\tau}$ is the submodule of $\Lambda^Q$ generated by
\[
\{tx+(1-t)y- (x \triangleright y),t^{-1}x + t^{-1} (t-1) y- (x \triangleright^{-1} y) \mid x,y \in Q\}.
\]
\end{enumerate}
\end{theorem}
\begin{proof}
This result follows from the discussion of \cite[Sec.\ 4]{mvaq1}, along with the right exactness of the functor $- \otimes_{\Lambda_{\mu}} \Lambda$. Alternatively, set all $t_i$ equal to $t$ in that discussion. \end{proof}

Part 2 of Theorem \ref{qpres} implies that once the function $\kappa_D:Q \to \{1, \dots, \mu\}$ is adjusted as in part 1, the map $\phi_{\tau}:M^{\textup{red}}_A(L) \to I_{\mu}$ will be determined by the fact that $\phi_{\tau}(\sigma_{\tau}(q))=(t_i-1) \otimes 1\thickspace \allowbreak \forall q \in \kappa_D^{-1}(\{i\})$. We deduce the ``forward'' direction of Theorem \ref{main1}: if $L$ and $L'$ are links and $f:Q^{\textup{red}}_A(L) \to Q^{\textup{red}}_A(L')$ is an isomorphism, then after adjusting component indices in $L$ and $L'$ so that $f$ maps the $\kappa_D^{-1}(i)$ orbit of $Q^{\textup{red}}_A(L)$ to the $\kappa_D^{-1}(i)$ orbit of $Q^{\textup{red}}_A(L')$ for each $i \in \{1. \dots, \mu \}$, $f$ will extend to an isomorphism $M^{\textup{red}}_A(L) \cong M^{\textup{red}}_A(L')$ that is compatible with the $\phi_{\tau}$ maps of $L$ and $L'$. 

The next two results give us the ``backward'' direction of Theorem \ref{main1}.

\begin{lemma}
\label{gen}
Suppose $0 \in W \subseteq \Lambda$. Then $W=\Lambda$ if and only if $W$ is closed under the following operations.
\begin{enumerate}
    \item $(w_1,w_2) \mapsto tw_1+(1-t)w_2$
    \item $(w_1,w_2) \mapsto 1+tw_1+(t-1)w_2$
    \item $(w_1,w_2) \mapsto t^{-1}w_1+(1-t^{-1})w_2$
    \item $(w_1,w_2) \mapsto -t^{-1}+t^{-1}w_1+(t^{-1}-1)w_2$
\end{enumerate}
\end{lemma}
\begin{proof}
If $W=\Lambda$ then $W$ is closed under all binary operations defined on $\Lambda$.

For the converse, suppose $W$ is closed under the four listed operations. Using $w_2=0$ in operations 1 and 3, we see that $W$ is closed under multiplication by $t^{\pm 1}$. Combining operation 2 with multiplication by $t^{-1}$, we see that $W$ contains the following elements: $(0,0) \mapsto 1$, $1 \mapsto t^{-1}$, $(t^{-1},0) \mapsto 2$, $2 \mapsto 2t^{-1}$, $(2t^{-1},0) \mapsto 3$, and so on. Combining operation 4 with multiplication by $t$, we see that $W$ contains the following elements: $(0,0) \mapsto -t^{-1}$, $-t^{-1} \mapsto -1$, $(-1,0) \mapsto -2t^{-1}$, $-2t^{-1} \mapsto -2$, $(-2,0) \mapsto -3t^{-1}$, $-3t^{-1} \mapsto -3$, and so on. We conclude that $\mathbb Z \subseteq W$.

Closure under multiplication by $t^{\pm 1}$ implies that $W$ contains every monomial $mt^n$ with $m,n \in \mathbb Z$. Now, suppose $\lambda \in \Lambda$ is not a monomial. Say $\lambda = n_1t^a+n_2t^{a+1} + \dots + n_kt^{a+k}$ for some $k \geq 1 \in \mathbb Z$ and some $n_1, \dots, n_k \in \mathbb Z$, with $n_1 \neq 0 \neq n_k$. Using induction on $k$, we may presume that $W$ contains $m_1t^b+m_2t^{b+1} + \dots + m_jt^{b+j}$ whenever $j<k$. Then $W$ contains both $w_1=n_1t^{a-1}+n_2t^{a} + \dots + n_{k-2}t^{a+k-3}+(n_{k-1}+n_{k})t^{a+k-2}$ and $w_2=-n_kt^{a+k-1}$. As $W$ is closed under operation 1, $W$ contains $tw_1+(1-t)w_2=\lambda$.
\end{proof}

\begin{proposition}
\label{mqprime}$
Q^{\textup{red}}_A(L)=(\phi_{\tau})^{-1}(\{(t_1-1) \otimes 1, \dots, (t_{\mu}-1) \otimes 1\}).$
\end{proposition}
\begin{proof}
Let $S=(\phi_{\tau})^{-1}(\{(t_1-1) \otimes 1, \dots, (t_{\mu}-1) \otimes 1\})$, and let $D$ be a diagram of $L$. Suppose $s \in S$ and $a \in A(D)$ have $\phi_{\tau}(s)=(t_i-1) \otimes 1$ and $\kappa_D(a)= j$. Then
\[
\phi_{\tau}(s \triangleright (\gamma_D(a) \otimes 1)) = \phi_{\tau}(ts + (1-t) \cdot (\gamma_D(a) \otimes 1))
\]
\[
=t \cdot ((t_i-1) \otimes 1) + (1-t) \cdot ((t_j-1) \otimes 1)= (t_i-1) \otimes t + (t_j-1) \otimes (1-t)
\]
\[
 =(t_j(t_i-1)) \otimes 1 + ((1-t_i)(t_j-1)) \otimes 1 = (t_i-1) \otimes 1 \textup{,}
\]
so $s \triangleright (\gamma_D(a) \otimes 1) \in S$. Clearly $\gamma_D(a) \otimes 1 \in S \thickspace \allowbreak \forall a \in A(D)$, so $S$ contains the subquandle of the standard Alexander quandle on $M^{\textup{red}}_A(L)$ generated by $\{ \gamma_D(a) \otimes 1 \mid a \in A(D) \}$. That is, $Q_A^{\textup{red}}(L) \subseteq S$.

Verifying the opposite inclusion is more difficult. Recall that $Q^{\textup{red}}_A(L)$ is a subquandle of the standard Alexander quandle on $M^{\textup{red}}_A(L)$, so $Q^{\textup{red}}_A(L)$ is closed under the operations $\triangleright,\triangleright^{-1}$ mentioned in Proposition \ref{standardq}. 

Suppose $a_1$ and $a_2$ are any two arcs of $D$. Let 
\[
W(a_1,a_2) = \{ w \in \Lambda \mid \gamma_D(a_1) \otimes 1 + w(t-1) \cdot ((\gamma_D(a_1)-\gamma_D(a_2)) \otimes 1) \in Q_A^{\textup{red}}(L) \} \textup{,}
\]
and let $W=W(a_1,a_2) \cap W(a_2,a_1)$. 

Note that if $w_1,w_2 \in W$ then $Q^{\textup{red}}_A(L)$ contains both
\[
(\gamma_D(a_1) \otimes 1 + w_1(t-1) \cdot ((\gamma_D(a_1)-\gamma_D(a_2)) \otimes 1))
\]
\[
 \triangleright (\gamma_D(a_1) \otimes 1 + w_2(t-1) \cdot ((\gamma_D(a_1)-\gamma_D(a_2)) \otimes 1))=
\]
\[
t \cdot (\gamma_D(a_1) \otimes 1 + w_1(t-1) \cdot ((\gamma_D(a_1)-\gamma_D(a_2)) \otimes 1))
\]
\[
+(1-t) \cdot (\gamma_D(a_1) \otimes 1 + w_2(t-1) \cdot ((\gamma_D(a_1)-\gamma_D(a_2)) \otimes 1))=
\]
\[
\gamma_D(a_1) \otimes 1 + (tw_1 + (1-t)w_2)(t-1) \cdot ((\gamma_D(a_1)-\gamma_D(a_2)) \otimes 1))
\]
and the element obtained from this by interchanging $a_1$ and $a_2$. It follows that $tw_1+(1-t)w_2 \in W$. That is, $W$ is closed under operation 1 of Lemma \ref{gen}.

For operation 2 of Lemma \ref{gen}, note that if $w_1,w_2 \in W$ then $Q^{\textup{red}}_A(L)$ contains both
\[
(\gamma_D(a_1) \otimes 1 + w_1(t-1) \cdot ((\gamma_D(a_1)-\gamma_D(a_2)) \otimes 1))
\]
\[
 \triangleright (\gamma_D(a_2) \otimes 1 + w_2(t-1) \cdot ((\gamma_D(a_2)-\gamma_D(a_1)) \otimes 1))=
\]
\[
t \cdot (\gamma_D(a_1) \otimes 1 + w_1(t-1) \cdot ((\gamma_D(a_1)-\gamma_D(a_2)) \otimes 1))
\]
\[
+(1-t) \cdot (\gamma_D(a_2) \otimes 1 + w_2(t-1) \cdot ((\gamma_D(a_2)-\gamma_D(a_1)) \otimes 1))=
\]
\[
\gamma_D(a_1) \otimes 1 + (t-1)(\gamma_D(a_1) \otimes 1)+tw_1(t-1) \cdot ((\gamma_D(a_1)-\gamma_D(a_2)) \otimes 1)
\]
\[
+(1-t) \cdot (\gamma_D(a_2) \otimes 1) + (1-t)w_2(t-1) \cdot ((\gamma_D(a_2)-\gamma_D(a_1)) \otimes 1))=
\]
\[
\gamma_D(a_1) \otimes 1 + (1+tw_1+(t-1)w_2)\cdot (t-1)((\gamma_D(a_1)-\gamma_D(a_2)) \otimes 1)
\]
and the element obtained from this by interchanging $a_1$ and $a_2$. It follows that $W$ is closed under the operation $(w_1,w_2) \mapsto 1+tw_1+(t-1)w_2$.

To show that $W$ is closed under operations 3 and 4 of Lemma \ref{gen}, use $\triangleright^{-1}$ instead of $\triangleright$:
\[
(\gamma_D(a_1) \otimes 1 + w_1(t-1) \cdot ((\gamma_D(a_1)-\gamma_D(a_2)) \otimes 1))
\]
\[
 \triangleright^{-1} (\gamma_D(a_1) \otimes 1 + w_2(t-1) \cdot ((\gamma_D(a_1)-\gamma_D(a_2)) \otimes 1))=
\]
\[
t^{-1} \cdot (\gamma_D(a_1) \otimes 1 + w_1(t-1) \cdot ((\gamma_D(a_1)-\gamma_D(a_2)) \otimes 1))
\]
\[
+(1-t^{-1}) \cdot (\gamma_D(a_1) \otimes 1 + w_2(t-1) \cdot ((\gamma_D(a_1)-\gamma_D(a_2)) \otimes 1))=
\]
\[
\gamma_D(a_1) \otimes 1 + (t^{-1}w_1 + (1-t^{-1})w_2)(t-1) \cdot ((\gamma_D(a_1)-\gamma_D(a_2)) \otimes 1)) \quad  \textup{and} 
\]
\[
(\gamma_D(a_1) \otimes 1 + w_1(t-1) \cdot ((\gamma_D(a_1)-\gamma_D(a_2)) \otimes 1))
\]
\[
 \triangleright^{-1} (\gamma_D(a_2) \otimes 1 + w_2(t-1) \cdot ((\gamma_D(a_2)-\gamma_D(a_1)) \otimes 1))=
\]
\[
t^{-1} \cdot (\gamma_D(a_1) \otimes 1 + w_1(t-1) \cdot ((\gamma_D(a_1)-\gamma_D(a_2)) \otimes 1))
\]
\[
+(1-t^{-1}) \cdot (\gamma_D(a_2) \otimes 1 + w_2(t-1) \cdot ((\gamma_D(a_2)-\gamma_D(a_1)) \otimes 1))=
\]
\[
\gamma_D(a_1) \otimes 1 + (-t^{-1}+t^{-1}w_1+(t^{-1}-1)w_2)\cdot (t-1)((\gamma_D(a_1)-\gamma_D(a_2)) \otimes 1).
\]

Lemma \ref{gen} now tells us that $W=\Lambda$. It follows that for every choice of $a_1,a_2 \in A(D)$ and $\lambda \in \Lambda$,
\[
q(a_1,a_2,\lambda)=\gamma_D(a_1) \otimes 1 + \lambda (t-1) \cdot ((\gamma_D(a_1)-\gamma_D(a_2)) \otimes 1) \in Q_A^{\textup{red}}(L).
\]

Now, suppose that $x \in Q_A^{\textup{red}}(L)$, $a_1,a_2 \in A(D)$ and $\lambda \in \Lambda$. Then $Q_A^{\textup{red}}(L)$ contains $x$, $q(a_1,a_2,\lambda)$ and $\gamma_D(a_2) \otimes 1$, so since $Q_A^{\textup{red}}(L)$ is closed under $\triangleright$ and $\triangleright^{-1}$, $Q_A^{\textup{red}}(L)$ also contains
\[
r(x,a_1,a_2,\lambda)=(x \triangleright^{-1} q(a_1,a_2,\lambda)) \triangleright (\gamma_D(a_2) \otimes 1) 
\]
\[
= t \cdot (x \triangleright^{-1} q(a_1,a_2,\lambda)) + (1-t) \cdot (\gamma_D(a_2) \otimes 1)
\]
\[
=t \cdot(t^{-1}x + (1-t^{-1})q(a_1,a_2,\lambda))  - (t-1) \cdot (\gamma_D(a_2) \otimes 1)
\]
\[
=x + (t-1) \cdot (q(a_1,a_2,\lambda) - (\gamma_D(a_2) \otimes 1) )
\]
\[
=x+ (t-1) \cdot ( 1 + \lambda (t-1) ) \cdot ((\gamma_D(a_1)-\gamma_D(a_2) \otimes 1)).
\]

Our next claim is that $x+ \lambda (t-1) \cdot ((\gamma_D(a_1)-\gamma_D(a_2)) \otimes 1) \in Q_A^{\textup{red}}(L)$ for all choices of $x \in Q_A^{\textup{red}}(L)$, $a_1,a_2 \in A(D)$ and $\lambda \in \Lambda$. If $\epsilon(\lambda)=1$, then $\lambda= 1 + \lambda'\cdot (t-1)$ for some $\lambda' \in \Lambda$, and 
\[
x+ \lambda (t-1) \cdot ((\gamma_D(a_1)-\gamma_D(a_2)) \otimes 1)=r(x,a_1,a_2,\lambda') \in Q_A^{\textup{red}}(L)\textup{,}
\]
so the claim is satisfied. If $\epsilon (\lambda) >1$ and
\[
x'=x+ (\lambda-1) (t-1) \cdot ((\gamma_D(a_1)-\gamma_D(a_2)) \otimes 1) \in Q_A^{\textup{red}}(L) \textup{,}
\]
then
\[
r(x',a_1,a_2,0)=x+ \lambda (t-1) \cdot ((\gamma_D(a_1)-\gamma_D(a_2)) \otimes 1) \in Q_A^{\textup{red}}(L) \textup{,}
\]
and the claim is satisfied. Using induction on $\epsilon(\lambda)$, we conclude that whenever $\epsilon(\lambda) \geq 1$, $x+ \lambda (t-1) \cdot ((\gamma_D(a_1)-\gamma_D(a_2)) \otimes 1) \in Q_A^{\textup{red}}(L)$. As $a_1$ and $a_2$ are arbitrary, it follows that the claim is also satisfied when $\epsilon(\lambda) \leq -1$, because 
\[
x+ \lambda (t-1) \cdot ((\gamma_D(a_1)-\gamma_D(a_2)) \otimes 1) = x+ (-\lambda) (t-1) \cdot ((\gamma_D(a_2)-\gamma_D(a_1)) \otimes 1).
\]

If $\epsilon(\lambda)=0$ then $\epsilon(\lambda-1)=-1$, so it follows that  $Q_A^{\textup{red}}(L)$ contains both $x' = x+ (\lambda-1) (t-1) \cdot ((\gamma_D(a_1)-\gamma_D(a_2)) \otimes 1) $ and $x+ \lambda (t-1) \cdot ((\gamma_D(a_1)-\gamma_D(a_2)) \otimes 1)$ $=r(x',a_1,a_2,0) $.
Again, the claim is satisfied.

According to Proposition \ref{kerphi}, this claim tells us that whenever $x \in Q_A^{\textup{red}}(L)$ and $y \in \ker \phi_{\tau}$, $x+y \in Q_A^{\textup{red}}(L)$. As $\phi_{\tau}(a) = (t_i-1) \otimes 1$ for every $a \in A(D)$ with $\kappa_D(a)=i$, it follows that $Q_A^{\textup{red}}(L)$ contains every element of the set $S$.
\end{proof}

Proposition \ref{mqprime} implies that if $L$ and $L'$ are $\phi_{\tau}$-equivalent links, then there is an isomorphism $M_A^{\textup{red}}(L) \cong M_A^{\textup{red}}(L')$ of $\Lambda$-modules, under which $Q_A^{\textup{red}}(L)$ and $Q_A^{\textup{red}}(L')$ correspond. This gives us the ``backward'' direction of Theorem \ref{main1}.

\section{Theorem \ref{maincor}}
\label{hopfex}

Suppose $L$ and $L'$ are links with $Q_A^{\textup{red}}(L) \cong Q_A^{\textup{red}}(L')$. Then Theorem \ref{main1} implies that after permuting the indices of components, $L$ and $L'$ will be $\phi_\tau$-equivalent. An isomorphism $f:\Mr(L) \to \Mr(L')$ with $\phi_\tau = \phi'_\tau \circ f$ will have $\fr = \eta \circ \phi_\tau = \eta \circ \phi'_\tau \circ f = \phi_{L'}^{\textup{red}} \circ f$, so $\ker \fr \cong \ker \phi_{L'}^{\textup{red}}$. We conclude that as a link invariant, $\Qr$ is at least as strong as the module $\ker \fr$. Examples verifying that $\Qr$ is strictly stronger than $\ker \fr$ are presented in the subsections below. 

Before discussing these examples, we prove that for knots, $\Qr(L)$ and $\ker \fr$ are equivalent invariants. We have already seen that $\Qr(L)$ determines $\ker \fr$. For the converse, recall that as noted in the introduction, every link $L$ has $\Mr(L) \cong \ker \fr \oplus \Lambda$. If we identify these two modules according to this isomorphism, then the map $\fr:\ker \phi_L^{\textup{red}} \oplus \Lambda \to \Lambda$ is given by a simple formula: $\fr(x,y) = y$. 

It follows that if $L$ and $L'$ are links and there is an isomorphism $f:\ker \fr \to \ker \phi_{L'}^{\textup{red}}$, then the resulting isomorphism $f \oplus \textup{id}:\ker \fr \oplus \Lambda \to \ker \phi_{L'}^{\textup{red}}\oplus \Lambda$ has $\fr = \phi_{L'}^{\textup{red}} \circ (f \oplus \textup{id})$. If $\mu=1$, then the map $\eta:I_\mu \otimes_{\Lambda_{\mu}} \Lambda \to \Lambda$ is an isomorphism, and 
\[
\phi_\tau = \eta^{-1} \circ \fr  = \eta^{-1} \circ \phi_{L'}^{\textup{red}} \circ (f \oplus \textup{id}) = \phi'_\tau \circ (f \oplus \textup{id})\text{,}
\]
so $L$ and $L'$ are $\phi_\tau$-equivalent. Theorem \ref{main1} then tells us that $\Qr(L) \cong \Qr(L')$.

To complete a proof of Theorem \ref{maincor}, it suffices to exhibit a pair of links $L,L'$ with $\ker \phi_{L}^{\textup{red}} \cong \ker \phi_{L'}^{\textup{red}}$ and $Q_A^{\textup{red}}(L) \not \cong Q_A^{\textup{red}}(L')$. In the rest of this section we present two such pairs.

\subsection{Two 4-component links}

\begin{figure} [bht]
\centering
\begin{tikzpicture} [>=angle 90]
\draw [thick] [->] (0,1.5) -- (-0.7,0.8);
\draw [thick] (-0.7,0.8) -- (-1.5,0);
\draw [thick] (0,-1.5) -- (1,-0.5);
\draw [thick] (0,1.5) -- (1.5,0);
\draw [thick] (1.5,0) -- (1.2,-0.3);
\draw [thick] (-1.5,0) -- (0,-1.5);
\draw [thick] (0.7,0) -- (2.2,-1.5);
\draw [thick] (0.7,0) -- (1,0.3);
\draw [thick] [->] (2.2,1.5) -- (1.5,0.8);
\draw [thick] (1.5,0.8) -- (1.2,0.5);
\draw [thick] (0.7,0) -- (2.2,-1.5);
\draw [thick] (2.2,1.5) -- (3.7,0);
\draw [thick] (3.7,0) -- (2.2,-1.5);
%right
\draw [thick] [->] (6.3+0,1.5) -- (6.3+-0.7,0.8);
\draw [thick] (6.3+-0.7,0.8) -- (6.3+-1.5,0);
\draw [thick] (6.3+0,-1.5) -- (6.3+1,-0.5);
\draw [thick] (6.3+0,1.5) -- (6.3+1.5,0);
\draw [thick] (6.3+1.5,0) -- (6.3+1.2,-0.3);
\draw [thick] (6.3+-1.5,0) -- (6.3+0,-1.5);
\draw [thick] (6.3+0.7,0) -- (6.3+2.2,-1.5);
\draw [thick] (6.3+0.7,0) -- (6.3+1,0.3);
\draw [thick] [->] (6.3+2.2,1.5) -- (6.3+1.5,0.8);
\draw [thick] (6.3+1.5,0.8) -- (6.3+1.2,0.5);
\draw [thick] (6.3+0.7,0) -- (6.3+2.2,-1.5);
\draw [thick] (6.3+2.2,1.5) -- (6.3+3.7,0);
\draw [thick] (6.3+3.7,0) -- (6.3+2.2,-1.5);

\node at (-1.1,0) {$a$};
\node at (3.3,0) {$b$};
\node at (5.2,0) {$c$};
\node at (9.6,0) {$d$};
\end{tikzpicture}
\caption{Two copies of the Hopf link.}
\label{hhfig}
\end{figure}
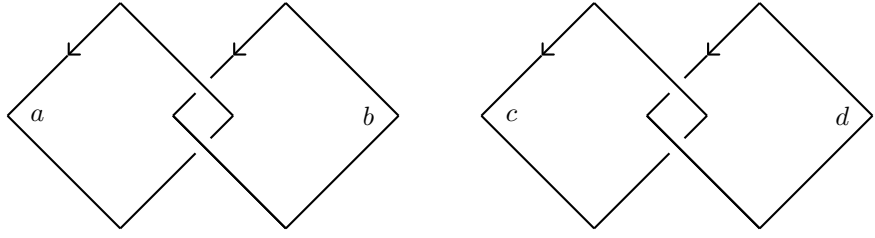

Let $L$ be the link illustrated in Fig.\ \ref{hhfig}. Then $M_A^{\textup{red}}(L)$ is generated by the four elements $\gamma_D(a) \otimes 1, \gamma_D(b) \otimes 1, \gamma_D(c) \otimes 1$ and $\gamma_D(d) \otimes 1$, subject to the crossing relations $(1-t) \cdot (\gamma_D(a) \otimes 1)=(1-t) \cdot (\gamma_D(b) \otimes 1)$ and $(1-t) \cdot (\gamma_D(c) \otimes 1)=(1-t) \cdot (\gamma_D(d) \otimes 1)$. It follows that 
\begin{equation}
\label{amods4}
M_A^{\textup{red}}(L) \cong \Lambda \oplus (\Lambda/(1-t)) \oplus
\Lambda \oplus (\Lambda/(1-t)) \textup,
\end{equation}
with the four summands generated by $\gamma_D(a) \otimes 1, (\gamma_D(b)-\gamma_D(a)) \otimes 1, (\gamma_D(c)-\gamma_D(a)) \otimes 1$ and $(\gamma_D(d)-\gamma_D(c)) \otimes 1$, respectively. The map $\fr:\Mr \to \Lambda$ maps all of $\gamma_D(a) \otimes 1,\gamma_D(b) \otimes 1,\gamma_D(c) \otimes 1,\gamma_D(d) \otimes 1$ to $1$, so 
\[
\ker \fr \cong  (\Lambda/(1-t)) \oplus
\Lambda \oplus (\Lambda/(1-t)).
\]

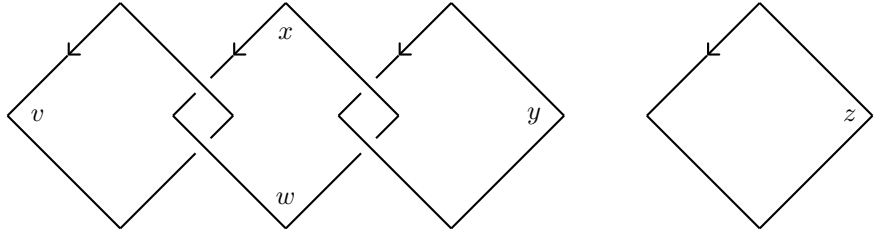
\begin{figure} [bht]
\centering
\begin{tikzpicture} [>=angle 90]
\draw [thick] [->] (0,1.5) -- (-0.7,0.8);
\draw [thick] (-0.7,0.8) -- (-1.5,0);
\draw [thick] (0,-1.5) -- (1,-0.5);
\draw [thick] (0,1.5) -- (1.5,0);
\draw [thick] (1.5,0) -- (1.2,-0.3);
\draw [thick] (-1.5,0) -- (0,-1.5);
\draw [thick] (0.7,0) -- (2.2,-1.5);
\draw [thick] (0.7,0) -- (1,0.3);
\draw [thick] [->] (2.2,1.5) -- (1.5,0.8);
\draw [thick] (1.5,0.8) -- (1.2,0.5);
\draw [thick] (2.2,1.5) -- (3.7,0);
\draw [thick] (3.7,0) -- (3.4,-0.3);
\draw [thick] (2.2,-1.5) -- (3.2,-0.5);
\draw [thick] (4.4,-1.5) -- (2.9,0);
\draw [thick] (4.4,-1.5) -- (5.9,0);
\draw [thick] (5.9,0) -- (4.4,1.5);
\draw [thick] [->] (4.4,1.5) -- (3.7,0.8);
\draw [thick] (3.7,0.8) -- (3.4,0.5);
\draw [thick] (3.2,0.3) -- (2.9,0);

\draw [thick] (8.5,-1.5) -- (7,0);
\draw [thick] (8.5,-1.5) -- (10,0);
\draw [thick] (8.5,1.5) -- (10,0);
\draw [thick] [->] (8.5,1.5) -- (7.8,0.8);
\draw [thick] (7,0) -- (7.8,0.8);
\node at (-1.1,0) {$v$};
\node at (2.2,-1.1) {$w$};
\node at (2.2,1.1) {$x$};
\node at (5.5,0) {$y$};
\node at (9.7,0) {$z$};
\end{tikzpicture}
\caption{This link is denoted $L'$ in Sec.\ 4.1.}
\label{hhufig}
\end{figure}

Now, let $L'$ be the link with the diagram $D'$ pictured in Fig.\ \ref{hhufig}. Then $M_A^{\textup{red}}(L')$ is generated by $\gamma_{D'}(v) \otimes 1, \gamma_{D'}(w) \otimes 1, \gamma_{D'}(x) \otimes 1,  \gamma_{D'}(y) \otimes 1$ and $\gamma_{D'}(z) \otimes 1$. The crossing relations from the two crossings on the left tell us that
\[
\gamma_{D'}(v) \otimes 1=(1-t) \cdot (\gamma_{D'}(w) \otimes 1)+ t\cdot (\gamma_{D'}(v) \otimes 1)
\]
\[
\textup{and} \quad \gamma_{D'}(w) \otimes 1=(1-t) \cdot (\gamma_{D'}(v) \otimes 1)+ t\cdot (\gamma_{D'}(x) \otimes 1) \textup,
\]
so $t\cdot (\gamma_{D'}(w) \otimes 1)=t\cdot (\gamma_{D'}(x) \otimes 1)$, and hence $\gamma_{D'}(w) \otimes 1=\gamma_{D'}(x) \otimes 1$. Taking this into account, the two relations from the crossings on the right are the same: $(1-t) \cdot (\gamma_{D'}(x) \otimes 1)=(1-t) \cdot (\gamma_{D'}(y) \otimes 1)$. It follows that
\begin{equation}
\label{amods5}
M_A^{\textup{red}}(L') \cong \Lambda \oplus (\Lambda/(1-t)) \oplus (\Lambda/(1-t)) \oplus \Lambda  \textup,
\end{equation}
with the four summands generated by $\gamma_{D'}(x) \otimes 1, (\gamma_{D'}(v)-\gamma_{D'}(x)) \otimes 1, (\gamma_{D'}(y)-\gamma_{D'}(x)) \otimes 1$ and $(\gamma_{D'}(z)-\gamma_{D'}(x)) \otimes 1$, respectively. Therefore
\[
\ker \phi_{L'}^{\textup{red}} \cong (\Lambda/(1-t)) \oplus (\Lambda/(1-t)) \oplus \Lambda .
\]

It is apparent that $\ker \fr \cong \ker \phi_{L'}^{\textup{red}}$. Nevertheless, we have the following.

\begin{proposition}
\label{nophitau}
No matter how their components are indexed, $L$ and $L'$ are not $\phi_\tau$-equivalent.
\end{proposition}
\begin{proof}
Recall that Lemma \ref{tenstruc} tells us 
\[
I_4 \otimes _{\Lambda_4} \Lambda \cong \Lambda \oplus \mathbb Z _ \epsilon \oplus \mathbb Z _ \epsilon \oplus \mathbb Z _ \epsilon \text{,}
\]
with the direct summands generated by $(t_1-1)\otimes 1, (t_2-t_1)\otimes 1, (t_3-t_1)\otimes 1$ and $(t_4-t_1)\otimes 1$, respectively. We use this isomorphism to identify $I_4 \otimes _{\Lambda_4} \Lambda$ with $\Lambda \oplus \mathbb Z _ \epsilon \oplus \mathbb Z _ \epsilon \oplus \mathbb Z _ \epsilon$.

Given an indexing of the components of $L$, the resulting map $\phi_{\tau}:\Mr(L) \to I_4 \otimes _{\Lambda_4} \Lambda$ will send $\gamma_D(a) \otimes 1, \gamma_D(b) \otimes 1, \gamma_D(c) \otimes 1$ and $\gamma_D(d) \otimes 1$ to $(1,0,0,0)$, $(1,1,0,0)$, $(1,0,1,0)$ and $(1,0,0,1)$, in some order. No matter what order is used, the image of $(\gamma_D(b)-\gamma_D(a)+\gamma_D(d)-\gamma_D(c)) \otimes 1$ will be of the form  $(0, \pm 1, \pm 1, \pm 1)$. Notice that according to (\ref{amods4}), $(\gamma_D(b)-\gamma_D(a)+\gamma_D(d)-\gamma_D(c)) \otimes 1$ is an element of $\Mr(L)$ that is annihilated by $t-1$.

Given an indexing of the components of $L'$, $\phi'_{\tau}:\Mr(L') \to I_4 \otimes _{\Lambda_4} \Lambda$ will map $\gamma_D(v) \otimes 1, \gamma_D(x) \otimes 1, \gamma_D(y) \otimes 1$ and $\gamma_D(z) \otimes 1$ to $(1,0,0,0)$, $(1,1,0,0)$, $(1,0,1,0)$ and $(1,0,0,1)$, in some order. According to (\ref{amods5}), every element of $\Mr(L')$ annihilated by $t-1$ is of the form 
\[
\lambda_1 (\gamma_D(v)-\gamma_D(x)) \otimes 1+ \lambda_2 (\gamma_D(y)-\gamma_D(x)) \otimes 1 
\]
\[
= (\lambda_1 \gamma_D(v)+ \lambda_2 \gamma_D(y)-(\lambda_1 + \lambda_2)\gamma_D(x)) \otimes 1
\]
for some $\lambda_1,\lambda_2 \in \Lambda$. It is easy to see that the image under $\phi'_{\tau}$ of such an element cannot be of the form $(0, \pm 1, \pm 1, \pm 1)$. 

Therefore, no isomorphism $f:\Mr(L) \to \Mr(L')$ has $\phi'_{\tau} \circ f = \phi_{\tau}$.
 \end{proof}

The proof of Proposition \ref{nophitau} is derived from an argument in the second paper of this series \cite[Sec.\ 15.4]{mvaq2}, where we observed that the links $L$ and $L'$ are distinguished by their involutory medial quandles. Results of \cite{mvaq2} imply that any two links with isomorphic reduced Alexander invariants and nonisomorphic involutory medial quandles must have $\mu \geq 4$. 

In the next subsection we see that in contrast, it is possible for medial quandles to distinguish 3-component links with isomorphic reduced Alexander invariants. We do not know of any analogous examples of 2-component links.

\subsection{Two 3-component links}

\begin{figure} [bth]
\centering
\begin{tikzpicture} [>=angle 90]
\draw [thick] [->] (0,0) -- (4.5,0);
\draw [thick] (4.5,0) -- (6,0);
\draw [thick] (6,0) -- (6,1);
\draw [thick] (6,1) -- (5.6,1.4);
\draw [thick] (5.4,1.6) -- (5,2);
\draw [thick] (5,1) -- (6,2);
\draw [thick] (6,2) -- (6,2.5);
\draw [thick] (5,1) -- (4.6,1.4);
\draw [thick] (4.4,1.6) -- (4,2);
\draw [thick] (4,1) -- (5,2);
\draw [thick] (3.6,2) -- (4,2);
\draw [thick] (3.6,1) -- (4,1);
\draw [thick] (3,2) -- (3.4,2);
\draw [thick] (3,1) -- (3.4,1);
\draw [thick] (3,2) -- (2,1);
\draw [thick] (3,1) -- (2.6,1.4);
\draw [thick] (2.4,1.6) -- (2,2);
\draw [thick] (1,2) -- (2,2);
\draw [thick] (1,1) -- (2,1);
\draw [thick] (0,1) -- (1,2);
\draw [thick] (0,2) -- (0.4,1.6);
\draw [thick] (1,1) -- (0.6,1.4);
\draw [thick] (0,2) -- (0,2.5);
\draw [thick] (0,0) -- (0,1);
\draw [thick] (0,2.5) -- (1.4,2.5);
\draw [thick] [->] (6,2.5) -- (4.5,2.5);
\draw [thick] (4.5,2.5) -- (1.6,2.5);
\draw [thick] [->] (1.5,3) -- (2.5,3);
\draw [thick] (3.5,3) -- (2.5,3);
\draw [thick] (3.5,3) -- (3.5,2.6);
\draw [thick] (3.5,2.4) -- (3.5,0.4);
\draw [thick] (1.5,0.4) -- (3.5,0.4);
\draw [thick] (1.5,0.4) -- (1.5,0.9);
\draw [thick] (1.5,1.1) -- (1.5,1.9);
\draw [thick] (1.5,2.1) -- (1.5,3);
\node at (5.8,2.7) {$a$};
\node at (0.2,2.8) {$b$};
\node at (2.25,1) {$c$};
\node at (3.9,1.8) {$d$};
\node at (1.3,1.5) {$g$};
\node at (1.3,2.9) {$e$};
\node at (2.5,0.6) {$f$};
\node at (5.8,0.2) {$j$};
\node at (3.9,1.2) {$h$};
\node at (3.15,1.2) {$i$};
\node at (1.3,0.8) {$*$};
\node at (2.8,1.5) {$*$};
\node at (4.5,1.2) {$*$};
\end{tikzpicture}
\caption{The link $L=L10n93\{0,1\}$ in the LinkInfo table \cite{linkinfo}.}
\label{thatlink}
\end{figure}
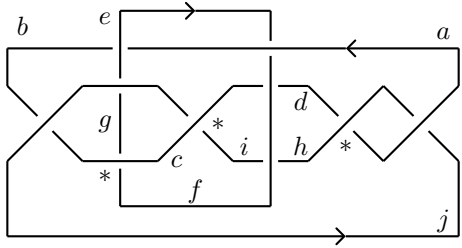

In this subsection, $L$ denotes the link pictured in Fig.\ \ref{thatlink}. We are grateful to Livingston, Moore and other researchers who developed and maintained the LinkInfo website \cite{linkinfo}, where we found this interesting example.

To describe $\Mr(L)$, we use the crossings not marked $*$ in Fig.\ \ref{thatlink} to eliminate generators other than $\overline a = \gamma_D(a) \otimes 1$, $\overline e = \gamma_D(e)\otimes 1$ and $\overline j = \gamma_D(j)\otimes 1$, as follows.
\begin{align*}
& \gamma_D(b)\otimes 1 = (1-t)\overline e + t\overline a \\
& \gamma_D(c)\otimes 1 = (1-t)\overline j + t(\gamma_D(b)\otimes 1) = t^2 \overline a + (t-t^2)\overline e + (1-t)\overline j\\
& \gamma_D(f)\otimes 1 = (1-t)\overline a + t\overline e\\
& \gamma_D(d)\otimes 1 = (1-t)(\gamma_D(f)\otimes 1) + t(\gamma_D(c)\otimes 1) \\
& \qquad \qquad \quad = (1-2t+t^2+t^3) \overline a + (t-t^3)\overline e + (t-t^2)\overline j\\
& \gamma_D(g)\otimes 1 = (1-t)\overline j + t\overline e\\
& \gamma_D(h)\otimes 1 = (1-t)\overline a + t\overline j\\
& \gamma_D(i)\otimes 1 = (1-t^{-1})(\gamma_D(f)\otimes 1) + t^{-1} (\gamma_D(h) \otimes 1) \\
& \qquad \qquad \quad = (1-t)\overline a  + (t-1)\overline e  + \overline j
\end{align*}

The three remaining crossings are marked $*$ in Fig.\ \ref{thatlink}. These three crossings provide three relations:
\begin{align*}
    & \gamma_D(g)\otimes 1=(1-t)(\gamma_D(c)\otimes 1)+t(\gamma_D(f)\otimes 1)\\
    & \overline j=(1-t)(\gamma_D(c)\otimes 1)+t(\gamma_D(i)\otimes 1) \\
    & \overline a = (1-t)(\gamma_D(h)\otimes 1)+t(\gamma_D(d)\otimes 1)
\end{align*}
The first two of these remaining crossing relations are both equivalent to
\begin{equation}
\label{rela}
0=(1-t) \cdot ((t^2+t)\overline a -t^2 \overline e -t \overline j) \text{,}
\end{equation}
and the third is equivalent to
\[
0=(t^2-1) \cdot ((t^2+t)\overline a -t^2 \overline e -t \overline j) \text{,}
\]
which is a consequence of (\ref{rela}). Therefore 
\begin{equation}
\label{isoone}
\Mr(L) \cong \Lambda \oplus \Lambda \oplus (\Lambda/(1-t)) \text{,}
\end{equation}
with the three summands generated by $\overline a$, $\overline e - \overline a$ and $\widehat j = (t^2+t)\overline a -t^2 \overline e -t \overline j$ (respectively). The map $\fr:\Mr \to \Lambda$ maps $\overline a$, $\overline e$ and $\overline j$ to $1$, so 
\[
\ker \fr \cong 
\Lambda \oplus (\Lambda/(1-t)).
\]

\begin{figure} [bht]
\centering
\begin{tikzpicture} [>=angle 90]
\draw [thick] [->] (0,1.5) -- (-0.7,0.8);
\draw [thick] (-0.7,0.8) -- (-1.5,0);
\draw [thick] (0,-1.5) -- (1,-0.5);
\draw [thick] (0,1.5) -- (1.5,0);
\draw [thick] (1.5,0) -- (1.2,-0.3);
\draw [thick] (-1.5,0) -- (0,-1.5);
\draw [thick] (0.7,0) -- (2.2,-1.5);
\draw [thick] (0.7,0) -- (1,0.3);
\draw [thick] [->] (2.2,1.5) -- (1.5,0.8);
\draw [thick] (1.5,0.8) -- (1.2,0.5);
\draw [thick] (0.7,0) -- (2.2,-1.5);
\draw [thick] (2.2,1.5) -- (3.7,0);
\draw [thick] (3.7,0) -- (2.2,-1.5);
%right
\draw [thick] [->] (6.3+0,1.5) -- (6.3+-0.7,0.8);
\draw [thick] (6.3+-0.7,0.8) -- (6.3+-1.5,0);
\draw [thick] (6.3+0,-1.5) -- (6.3+1.5,0);
\draw [thick] (6.3+0,1.5) -- (6.3+1.5,0);
%draw [thick] (6.3+1.5,0) -- (6.3+1,-0.5);
\draw [thick] (6.3+-1.5,0) -- (6.3+0,-1.5);

\node at (-1.1,0) {$x$};
\node at (3.3,0) {$y$};
\node at (5.2,0) {$z$};
\end{tikzpicture}
\caption{This link is denoted $L'$ in Sec.\ 4.2.}
\label{hhhufig}
\end{figure}
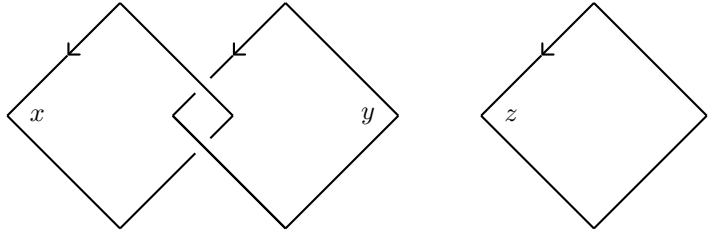

Now, let $L'$ be the link with the diagram $D'$ pictured in Fig.\ \ref{hhhufig}. Then it is easy to see that
\begin{equation}
\label{isotwo}
\Mr(L') \cong \Lambda \oplus (\Lambda/(1-t)) \oplus \Lambda \text{,}
\end{equation}
with the three summands generated by $\gamma_{D'}(x) \otimes 1$, $(\gamma_{D'}(y)-\gamma_{D'}(x)) \otimes 1$ and $(\gamma_{D'}(z)-\gamma_{D'}(x)) \otimes 1$ (respectively). As $\phi_{L'}^{\textup{red}}:\Mr(L') \to \Lambda$ maps $\gamma_{D'}(x) \otimes 1$, $\gamma_{D'}(y) \otimes 1$ and $\gamma_{D'}(z) \otimes 1$ to $1$, it follows that
\[
\ker \phi_{L'}^{\textup{red}} \cong 
\Lambda \oplus (\Lambda/(1-t)).
\]

The relationship between the links denoted $L$ and $L'$ in this subsection is similar to the relationship between the links denoted $L$ and $L'$ in Sec.\ 4.1.

\begin{proposition}
\label{nophitau2}
No matter how their components are indexed, $L$ and $L'$ are not $\phi_\tau$-equivalent.
\end{proposition}
\begin{proof}
Identify $I_3 \otimes _{\Lambda_3} \Lambda$ with $\Lambda_3 \oplus \mathbb Z _ \epsilon \oplus \mathbb Z _ \epsilon$ using the isomorphism of Lemma \ref{tenstruc}. Given an indexing of the components of $L$, the resulting map $\phi_{\tau}:\Mr(L) \to I_3 \otimes _{\Lambda_3} \Lambda$ will send $\overline a, \overline e$ and $\overline j$ to $(1,0,0)$, $(1,1,0)$, and $(1,0,1)$, in some order. 

According to (\ref{isoone}), every element $m \in \Mr(L)$ with $(1-t)m=0$ is of the form $m=\lambda \widehat j = \lambda ((t^2+t)\overline a -t^2 \overline e -t \overline j)$ for some $\lambda \in \Lambda$. If $\phi_ \tau(\overline a ) = (1,0,0)$, then $\phi_\tau(\lambda \widehat j) = \epsilon(\lambda)\phi_\tau((t^2+t)\overline a -t^2 \overline e -t \overline j) = (0,-\epsilon(\lambda),-\epsilon(\lambda))$. If $\phi_ \tau(\overline a ) = (0,1,0)$, then $\phi_\tau(\lambda \widehat j) = \epsilon(\lambda)\phi_\tau((t^2+t)\overline a -t^2 \overline e -t \overline j) = (0,2\epsilon(\lambda),-\epsilon(\lambda))$. If $\phi_ \tau(\overline a ) = (0,0,1)$, then $\phi_\tau(\lambda \widehat j) = \epsilon(\lambda)\phi_\tau((t^2+t)\overline a -t^2 \overline e -t \overline j) = (0,-\epsilon(\lambda),2\epsilon(\lambda))$. No matter how the components of $L$ are indexed, if $m \in \Mr(L)$ has $(1-t)m=0$, then $\phi_\tau(m) \neq (0,1,0)$.

Suppose the components of $L'=K'_1 \cup K'_2 \cup K'_3$ are indexed so that $x,y$ and $z$ are the images in $D'$ of $K'_1,K'_2$ and $K'_3$, in order. Then $\phi'_{\tau}:\Mr(L') \to I_3 \otimes _{\Lambda_3} \Lambda$ has $\phi'_{\tau}(\gamma_{D'}(x) \otimes 1) = (1,0,0)$, $\phi'_{\tau}(\gamma_{D'}(y) \otimes 1)=(1,1,0)$ and $\phi'_{\tau}(\gamma_{D'}(z) \otimes 1)=(1,0,1)$. It follows that $\phi'_{\tau}((\gamma_{D'}(y)-\gamma_{D'}(x)) \otimes 1)=(0,1,0)$. 

According to (\ref{isotwo}), $m=(\gamma_{D'}(y)-\gamma_{D'}(x)) \otimes 1$ is an element of $\Mr(L')$ with $(1-t)m=0$. Every isomorphism $f:\Mr(L) \to \Mr(L')$ maps elements of $\Mr(L)$ annihilated by $1-t$ to elements of $\Mr(L')$ annihilated by $1-t$, so every isomorphism $f:\Mr(L) \to \Mr(L')$ has $\phi_{\tau} \neq \phi'_{\tau} \circ f$. \end{proof}

All of the links discussed in this section illustrate a fact mentioned in the introduction: Theorem \ref{main1} does not hold without allowing for re-indexing of link components. The reasoning is simple. For instance, suppose the components of the link denoted $L'$ in this subsection are $K'_1, K'_2,K'_3$ in $x,y,z$ order, as in the proof of Proposition \ref{nophitau2}. The isomorphism (\ref{isotwo}) implies that if $m \in \Mr(L')$ has $(1-t)m=0$, then $m = \lambda \cdot ((\gamma_{D'}(y)-\gamma_{D'}(x)) \otimes 1)$ for some $\lambda \in \Lambda$, so $\phi'_\tau:\Mr(L') \to I_3 \otimes _{\Lambda_3} \Lambda$ has $\phi'_\tau(m)=(0,\epsilon(\lambda),0)$. Now, let $L''=K''_1 \cup K''_2 \cup K''_3$ be the same link, but with $K''_1, K''_2,K''_3$ in $x,z,y$ order. Then the same analysis leads to the conclusion that if $m \in \Mr(L'')$ has $(1-t)m=0$, the map $\phi''_\tau:\Mr(L'') \to I_3 \otimes _{\Lambda_3} \Lambda$ has $\phi''_\tau(m)=(0,0,n)$ for some integer $n$. It follows that no isomorphism between the modules $\Mr(L')$ and $\Mr(L'')$ is compatible with $\phi'_\tau$ and $\phi''_\tau$, so $L'$ and $L''$ are not $\phi_\tau$-equivalent. However $\Qr(L') \cong \Qr(L'')$, because the only difference between $L'$ and $L''$ is the indexing of their components, and the quandles do not detect the indices.

\section{Medial Quandles}
\label{struc}

In this section, we summarize some ideas from the general theory of quandles. We refer to the work of Jedli\v{c}ka, Pilitowska, Stanovsk\'{y} and Zamojska-Dzienio \cite{JPSZ1, JPSZ2} for a more thorough discussion. Proofs are included for the reader's convenience.

If $Q$ is a quandle and $y \in Q$, then the map $\beta_y$ is a \emph{translation} of $Q$. Part 2 of Definition \ref{qdef} implies that $\beta_y$ has an inverse function; the notation $\beta_y^{-1}(x) = x \triangleright ^{-1} y$ is often used. Notice that part 3 of Definition \ref{qdef} can be written as $\beta_z(x \triangleright y)= \beta_z(x) \triangleright \beta_z(y)$, so every translation of $Q$ is a quandle automorphism. (Indeed, some authors call translations ``inner automorphisms.'') Of course it follows that the inverse function of a translation is an automorphism too, so $\beta_z^{-1}(x \triangleright^{-1} y)= \beta_z^{-1}(x) \triangleright^{-1} \beta_z(y) \thickspace \allowbreak \forall x,y,z \in Q$. This implies that $\triangleright^{-1}$ also defines a quandle structure on the set $Q$. 

The fact that $\beta_z$ is an automorphism of $Q$ also implies that $\beta_z(x \triangleright^{-1} y)=\beta_z(x) \triangleright^{-1} \beta_z(y) \thickspace \allowbreak \forall x,y,z \in Q$. That is, $\beta_z \beta_y^{-1} (x)=\beta_{\beta_z(y)}^{-1} \beta_z(x) \thickspace \allowbreak \forall x,y,z \in Q$. It follows that there is a special way to express conjugation of translations in $\Aut(Q)$: $\beta_z \beta_y \beta_z^{-1}=\beta_{\beta_z(y)}$.

If $y,z \in Q$ then the composition $\beta_y\beta_z^{-1}$ is an \emph{elementary displacement} of $Q$. The subgroup of the automorphism group $\Aut(Q)$ generated by the elementary displacements is denoted $\Dis(Q)$, and its elements are called \emph{displacements}. (Some references use the term ``transvections'' instead.) 

\begin{proposition} 
\label{disprop}
Here are six properties of displacements and translations.
\begin{enumerate}
    \item If $d \in \Dis(Q)$ then $d=\beta_{y_1}\beta^{-1}_{y_2} \dots \beta_{y_{2n-1}}\beta^{-1}_{y_{2n}}$ for some $y_1,\dots,y_{2n} \in Q$.
    \item The elementary displacements also include products of the form $\beta_y^{-1} \beta_z$.
    \item If $y_1,\dots,y_n \in Q$, $m_1, \dots, m_n \in \mathbb Z$ and $\sum m_i = 0$, then $\prod \beta_{y_i}^{m_i} \in \Dis(Q)$.
    \item If $f \in \Aut(Q)$ and $y \in Q$, then $\beta_{f(y)} = f\beta_y f^{-1}$.
    \item $\Dis(Q)$ is a normal subgroup of $\Aut(Q)$.
    \item If $f:Q_1 \to Q_2$ is a surjective quandle map, then $f$ induces a surjective homomorphism $\Dis(f):\Dis(Q_1) \to \Dis(Q_2)$, defined by 
    \[
    \Dis(f)(\prod \beta_{y_i}^{m_i}) =  \prod \beta_{f(y_i)}^{m_i}
    \]
    whenever $y_1,\dots,y_n \in Q_1$, $m_1, \dots, m_n \in \mathbb Z$ and $\sum m_i = 0$.
\end{enumerate}
\end{proposition}
\begin{proof}
For the first property, note that the inverse of an elementary displacement is an elementary displacement. It follows that every element of $Dis(Q)$ is a product of elementary displacements.

For the second property, replace $y$ with $y'=\beta_z^{-1}(y)$ in the equality $\beta_z \beta_y^{-1} =\beta_{\beta_z(y)}^{-1} \beta_z$ mentioned above. The result is $\beta_z \beta^{-1}_{y'} = \beta_y^{-1} \beta_z$. 

For the third property, introduce repetitions in the list $y_1, \dots, y_n$ so that $m_1, \dots, m_n \in \{-1,1\}$. If $m_1 \neq m_2$ then $\prod \beta_{y_i}^{m_i}=(\beta_{y_1}^{m_1}\beta_{y_2}^{m_2})(\prod_{i>2} \beta_{y_i}^{m_i})$, and induction on $n$ applies. If $m_1 = m_2$ then find the least index $j$ with $m_1 \neq m_j$, and apply the second property to replace $\beta_{y_{j-1}}^{m_{j-1}}$ and $\beta_{y_j}^{m_j}$ with  $\beta_{y'_{j-1}}^{m_j}$ and $\beta_{y'_j}^{m_{j-1}}$, respectively, so that $\beta_{y'_{j-1}}^{m_j}\beta_{y'_j}^{m_{j-1}} = \beta_{y_{j-1}}^{m_{j-1}}\beta_{y_j}^{m_j}$. Now induction on $j$ applies.

For the fourth property, notice that for each $x \in Q$, $f\beta_y(x) = f(x \triangleright y) = f(x) \triangleright f(y) = \beta_{f(y)} f (x)$. Hence $f\beta_y = \beta_{f(y)} f$. 

For the fifth property, notice that if $\beta_y \beta_z^{-1}$ is an elementary displacement, then for every $f \in \Aut(Q)$,
\[
f \beta_y \beta_z^{-1} f^{-1}=f \beta_y (f \beta_z)^{-1} = \beta_{f(y)} f (\beta_{f(z)} f)^{-1} = \beta_{f(y)} f f^{-1} \beta_{f(z)}^{-1}= \beta_{f(y)} \beta_{f(z)}^{-1}
\]
is also an elementary displacement. 

For the sixth property, note that as $f$ is a quandle map, $f(x \triangleright y) = f(x) \triangleright f(y)\thickspace \allowbreak \forall x,y \in Q_1$. Also, if $x,y \in Q_1$ then $f(x \triangleright ^{-1} y) \triangleright f(y) = f((x \triangleright^{-1} y) \triangleright y)=f(x)$, so $f(x \triangleright ^{-1} y) = f(x) \triangleright^{-1} f(y)$. We deduce that $f \circ \beta_y^{\pm 1} = \beta_{f(y)}^{\pm 1} \circ f$  for every $y \in Q_1$. This implies that $f \circ \beta_y^n = \beta_{f(y)}^n \circ f\thickspace \allowbreak \forall y \in Q_1\thickspace \allowbreak \forall n \in \mathbb Z$.

Suppose $y_1, \dots, y_n, z_1, \dots, z_p \in Q_1$, $m_1,\dots,m_n,\ell_1,\dots, \ell_p \in \mathbb Z$, and 
\[
    \prod_{i=1}^n \beta_{y_i}^{m_i} =  \prod_{j=1}^p \beta_{z_j}^{\ell_j}.
\]
Then 
\[
     \left( \prod_{i=1}^n \beta_{f(y_i)}^{m_i} \right) \circ f = \left( \prod_{i=1}^{n-1} \beta_{f(y_i)}^{m_i} \right) \circ f \circ \beta_{y_n}^{m_n} = \dots = f \circ \left( \prod_{i=1}^n \beta_{y_i}^{m_i} \right)
    \]
    \[
     = f \circ \left( \prod_{j=1}^p \beta_{z_j}^{\ell_j} \right) = \beta_{f(z_1)}^{\ell_1} \circ f \circ \left( \prod_{j=2}^p \beta_{z_j}^{\ell_j} \right) = \dots = \left( \prod_{j=1}^p \beta_{f(z_j)}^{\ell_j} \right) \circ f.
    \]
As $f$ is surjective, we deduce that 
\[
 \prod_{i=1}^n \beta_{y_i}^{m_i} =  \prod_{j=1}^p \beta_{z_j}^{\ell_j}
 \implies \prod_{i=1}^n \beta_{f(y_i)}^{m_i} = \prod_{j=1}^p \beta_{f(z_j)}^{\ell_j}.
\]
That is, we have a well-defined function mapping products of powers of translations of $Q_1$ to products of powers of translations of $Q_2$, according to the formula 
\[
\prod_{i=1}^n \beta_{y_i}^{m_i} \mapsto \prod_{i=1}^n \beta_{f(y_i)}^{m_i}.
\]
The restriction of this function to products with $\sum m_i = 0$ is the map $\Dis(f)$ mentioned in the sixth property of the proposition. As the function maps products to products, it is obvious that it is a homomorphism. The image of $\Dis(f)$ includes every elementary displacement of $Q_2$, so $\Dis(f)$ is surjective. 
\end{proof}

If $x$ is an element of a quandle $Q$, then the orbit of $x$ in $Q$ is the smallest subset that contains $x$ and is preserved by $\beta_y$ and $\beta^{-1}_y$, for every $y \in Q$.

\begin{proposition}
\label{orb}
If $x \in Q$, then the orbit of $x$ in $Q$ is $\{d(x) \mid d \in \Dis(Q)\}$.
\end{proposition}
\begin{proof}
A displacement is a composition of translations and their inverses, so the orbit of $x$ includes $d(x)$ for every displacement $d$. 

Now, suppose $y$ is an element of the orbit of $x$. Then there are $y_1,\dots,y_n \in Q$ and $\epsilon _1, \dots, \epsilon _n \in \{-1,1\}$ such that $y=\beta^{\epsilon _n}_{y_n}\cdots \beta^{\epsilon _1}_{y_1}(x)$. For $1 \leq i \leq n$, let $x_i=\beta^{\epsilon _i}_{y_i}\cdots \beta^{\epsilon _1}_{y_1}(x)$. Then $\beta_{x_i}(x_i)=\beta^{-1}_{x_i}(x_i)=x_i$ for every $i$, so the function
\[
d=(\beta^{-\epsilon _n}_{x_n}\beta^{\epsilon _n}_{y_n}) (\beta^{-\epsilon _{n-1}}_{x_{n-1}}\beta^{\epsilon _{n-1}}_{y_{n-1}})\cdots (\beta^{-\epsilon _1}_{x_1}\beta^{\epsilon _1}_{y_1})
\]
has $y=d(x)$. Each product $\beta^{-\epsilon _i}_{x_i}\beta^{\epsilon _i}_{y_i}$ is an elementary displacement, so $d \in \Dis(Q)$. 
\end{proof}

\begin{definition}
\label{semireg}
A quandle is \emph{semiregular} if the only displacement with a fixed point is the identity map.
\end{definition}

According to Proposition \ref{orb}, if $Q$ is semiregular then for each $x \in Q$, the map $d \mapsto d(x)$ is a bijection from $\Dis(Q)$ to the orbit of $x$ in $Q$.

\begin{definition}
A quandle is \emph{medial} if it has the property that $(w\triangleright x) \triangleright (y\triangleright z)=(w\triangleright y) \triangleright (x\triangleright z)\thickspace \allowbreak \forall w,x,y,z \in Q$.
\end{definition}

\begin{proposition}
\label{medprop}
These three properties of a quandle $Q$ are equivalent.
\begin{enumerate}
    \item $Q$ is medial. 
    \item $\beta_q \beta^{-1}_r \beta_s=\beta_s \beta^{-1}_r \beta_q \thickspace \allowbreak \forall q,r,s \in Q$.
    \item $\Dis(Q)$ is an abelian group.
\end{enumerate}
\end{proposition}
\begin{proof} Suppose first that $Q$ is medial. Recall that as mentioned in the third paragraph of this section, if $q,r \in Q$ then the formula $\beta_{q \triangleright r}=\beta_{\beta_ r(q)}= \beta_r\beta_q\beta_r^{-1}$ holds. Also, if $q,r,s \in Q$ then for every $x \in Q$,
\[
\beta_{q \triangleright r}(\beta_s(x))= (x \triangleright s) \triangleright (q \triangleright r)= (x \triangleright q) \triangleright (s \triangleright r) = \beta_{s \triangleright r}(\beta_q(x)) \textup{,}
\]
so $\beta_{q \triangleright r}\beta_s=\beta_{s \triangleright r}\beta_q$.  Using the first formula twice, we obtain the second property:
\[
\beta_q \beta^{-1}_r \beta_s = \beta_r ^{-1} \cdot \beta_r \beta_q \beta^{-1}_r \beta_s =  \beta_r ^{-1} \cdot  \beta_{q \triangleright r}\beta_s 
\]
\[
 =  \beta_r ^{-1} \cdot \beta_{s \triangleright r}\beta_q =  \beta_r ^{-1} \cdot  \beta_r\beta_s\beta_r^{-1}\beta_q   =  \beta_s\beta_r^{-1}\beta_q.
\]

Now, suppose the second property holds. Let $\beta_q \beta_r^{-1}$ and $\beta_s \beta_t^{-1}$ be elementary displacements of $Q$. Then using the second property twice, 
\[
(\beta_q \beta_r^{-1})(\beta_s \beta_t^{-1}) = (\beta_q \beta_r^{-1}\beta_s) \beta_t^{-1}  = (\beta_s \beta_r^{-1}\beta_q) \beta_t^{-1}   = \beta_s (\beta_r^{-1}\beta_q \beta_t^{-1})
\]
\[
= \beta_s (\beta_t\beta_q^{-1} \beta_r)^{-1} = \beta_s (\beta_r\beta_q^{-1} \beta_t)^{-1}  = \beta_s (\beta_t^{-1}\beta_q \beta_r^{-1})  = (\beta_s \beta_t^{-1})(\beta_q \beta_r^{-1}).
\]
That is, the elementary displacements $\beta_q \beta_r^{-1}$ and $\beta_s \beta_t^{-1}$ commute with each other. As the group $\Dis(Q)$ is generated by elementary displacements, it follows that $\Dis(Q)$ is abelian.

For the last part of the proof, suppose $\Dis(Q)$ is abelian. Let $a,b,c,d \in Q$, and let $x=\beta_d^{-1} \beta_b(a)$. Then using part 3 of Definition \ref{qdef} twice,
\[
(a \triangleright b ) \triangleright (c \triangleright d) = (x \triangleright d) \triangleright (c \triangleright d) = (x  \triangleright c) \triangleright d = \beta_d \beta_c \beta_d^{-1} \beta_b(a)
\]
\[
= \beta_d (\beta_c \beta_d^{-1}) (\beta_b \beta_d^{-1}) \beta_d (a) = \beta_d (\beta_b \beta_d^{-1}) (\beta_c \beta_d^{-1}) \beta_d (a) = \beta_d \beta_b \beta_d^{-1} \beta_c (a)
\]
\[
= (\beta_d^{-1} (a \triangleright c) \triangleright b) \triangleright d = (\beta_d^{-1} (a \triangleright c) \triangleright d) \triangleright (b \triangleright d) = (a \triangleright c) \triangleright (b \triangleright d) \textup{,}
\]
so $Q$ is medial.
\end{proof}

\begin{corollary}
\label{dismod}
Let $Q$ be a medial quandle, and let $q^*$ be a fixed element of $Q$. Then  $\Dis(Q)$ is a $\Lambda$-module, with scalar multiplication given by 
    \[
    t \cdot d = \beta_{q^*} \, d \, \beta^{-1}_{q^*}  \thickspace \allowbreak \forall d \in \Dis(Q).
    \]
Changing the choice of $q^*$ does not change the $\Lambda$-module structure of $\Dis(Q)$.
\end{corollary}
\begin{proof}
If $A$ is an abelian group with an automorphism $\alpha:A \to A$, we obtain a $\Lambda$-module structure on $A$ by setting $t \cdot a = \alpha(a)\thickspace \allowbreak \forall a \in A$. The first assertion follows, because conjugation by $\beta_{q^*}$ is an automorphism of $\Aut(Q)$, and it defines an automorphism of the normal abelian subgroup $\Dis(Q) \subset \Aut(Q)$. 

Suppose $q^*$ and $q^{**}$ are fixed elements of $Q$. As $\Dis(Q)$ is abelian, we have
\[
\beta_{q^*} \, d \, \beta^{-1}_{q^*} = \beta_{q^*} (\beta^{-1}_{q^*}\beta_{q^{**}})(\beta^{-1}_{q^{**}}\beta_{q^*}) \, d \, \beta^{-1}_{q^*} = \beta_{q^*} (\beta^{-1}_{q^*}\beta_{q^{**}}) \, d \,  (\beta^{-1}_{q^{**}}\beta_{q^*}) \beta^{-1}_{q^*} 
\]
\[
= \beta_{q^{**}} \, d \, \beta^{-1}_{q^{**}} 
\]
for every $d \in \Dis(Q)$.
\end{proof}

In the proof of the sixth part of Proposition \ref{disprop}, we showed that if $f:Q_1 \to Q_2$ is a surjective quandle map, then $f$ induces a well-defined function mapping products of translations of $Q_1$ to products of translations of $Q_2$, given by 
\[
\prod_{i=1}^n \beta_{y_i}^{m_i} \mapsto \prod_{i=1}^n \beta_{f(y_i)}^{m_i}.
\]
If $Q_1$ and $Q_2$ are medial quandles, then the induced function $\Dis(f):\Dis(Q_1) \to \Dis(Q_2)$ is not only a homomorphism of abelian groups. Corollary \ref{dismod} defines scalar multiplication in $\Dis(Q_1)$ and $\Dis(Q_2)$ using multiplication of translations, so $\Dis(f)$ is actually a homomorphism of $\Lambda$-modules.

Recall that as discussed in Sec.\ \ref{proof1}, if $M$ is a $\Lambda$-module then the \emph{standard Alexander quandle} on $M$ is given by the operation $x \triangleright y = tx + (1-t)y$.

\begin{proposition}
\label{coreprop}
If $M$ is a $\Lambda$-module, then the standard Alexander quandle on $M$ is a semiregular medial quandle. Also, $\Dis(M) \cong (1-t)M$ as $\Lambda$-modules.
\end{proposition}
\begin{proof}
 It is easy to see that $M$ is a quandle, with $\triangleright ^{-1}$ given by $x \triangleright^{-1} y = t^{-1}x + (1-t^{-1})y$. The medial property is verified at the beginning of Sec.\ \ref{proof1}.
 
 Notice that if $y,z \in M$ then for every $x \in M$, 
 \[
 \beta_y \beta_z^{-1}(x)=\beta_y(t^{-1}x+(1-t^{-1})z)
 \]
 \[
 =t(t^{-1}x+(1-t^{-1})z)+(1-t)y=x+(1-t)(y-z).
 \]
It follows that there is a well-defined function $g:(1-t)M \to \Dis(M)$, with $g((1-t)m)$ being the displacement given by $g((1-t)m)(x) = x+ (1-t)m$. (That is, $g((1-t)m)=\beta_m\beta^{-1}_0$.) It is obvious that $g$ is injective, as $(1-t)m = g((1-t)m)(0)$. 

If $m_1,m_2 \in M$ then $g((1-t)m_1+(1-t)m_2)$ is the function with
\[
g((1-t)m_1+(1-t)m_2)(x) = x + (1-t)m_1+(1-t)m_2
\]
\[
= (x + (1-t)m_2)+(1-t)m_1 = g((1-t)m_1)(g((1-t)m_2)(x)) \textup{,}
\]
so $g((1-t)m_1+(1-t)m_2) = g((1-t)m_1) \circ g((1-t)m_2)$. That is, $g$ is a homomorphism of abelian groups. 

Moreover, the image of $g$ contains every elementary displacement, because $\beta_y \beta_z^{-1}$ $=g((1-t)(y-z))$. It follows that $g$ is surjective.

Choose a fixed element $q^* \in M$, and define a $\Lambda$-module structure on $\Dis(M)$ using $t \cdot d = \beta_{q^*} \, d \, \beta^{-1}_{q^*}$, as in Corollary \ref{dismod}. Then for all $m,x \in M$,
\[
g(t \cdot (1-t)m)(x) = t(1-t)m + x = x+(t-1)q^*+t(1-t)m + (1-t)q^*
\]
\[
= (t^{-1} \cdot (x+(t-1)q^*+t(1-t)m)) \triangleright q^*
=\beta_{q^*}(t^{-1} \cdot (x+(t-1)q^*)+(1-t)m)
\]
\[
=(\beta_{q^*} \circ g((1-t)m))(t^{-1} \cdot (x+(t-1)q^*))= (\beta_{q^*} \circ g((1-t)m))(x \triangleright^{-1} q^*)
\]
\[
= (\beta_{q^*} \circ g((1-t)m) \circ \beta_{q^*}^{-1})(x) = (t \cdot g((1-t)m))(x).
\]
Hence $g(t \cdot (1-t)m)=t \cdot g( (1-t)m)$, so $g$ is an isomorphism of $\Lambda$-modules.

To verify semiregularity of $M$, suppose $d \in \Dis(M)$. Then $d=g((1-t)m)$ for some $m \in M$. If there is an $x \in M$ with $x=d(x)$, then $x = g((1-t)m)(x) = x+ (1-t)m$, so $(1-t)m=0$. Therefore $d=g(0)$ is the identity map of $M$. \end{proof}

Notice that if $M$ is a standard Alexander quandle then the surjectivity of the map $g$ used in the proof of Proposition  \ref{coreprop} implies that every displacement of $M$ is of the form $g((1-t)m) = \beta_m \beta_0^{-1}$ for some $m \in M$. In the terminology of Jedli\v{c}ka, Pilitowska, Stanovsk\'{y} and Zamojska-Dzienio \cite{JPSZ2}, standard Alexander modules have ``tiny'' displacement groups.

Standard Alexander quandles are the building blocks of medial quandles.

\begin{proposition}
\label{orbstruc}
Let $x$ be an element of a medial quandle $Q$, and let $Q_x$ be the orbit of $x$ in $Q$. Then $\Fix(x)=\{d \in \Dis(Q) \mid d(x)=x\}$ is a $\Lambda$-submodule of $\Dis(Q)$, and $Q_x$ is isomorphic, as a quandle, to the standard Alexander quandle on the quotient module $\Dis(Q)/\Fix(x)$. A quandle isomorphism $e:\Dis(Q)/\Fix(x) \to Q_x$ is given by $e(d+\Fix(x)) = d(x)$.
\end{proposition}
\begin{proof}
It is easy to see that $\Fix(x)$ is a subgroup of $\Dis(Q)$. According to Corollary \ref{dismod}, if $d \in \Fix(x)$ then $(t \cdot d)(x) = \beta_x d\beta_x^{-1}(x) = \beta_x d(x) = \beta_x (x) = x$ and $(t^{-1} \cdot d)(x) = \beta_x^{-1} d\beta_x(x) = \beta_x^{-1} d(x) = \beta_x^{-1} (x) = x$. Therefore $\Fix(x)$ is closed under scalar multiplication, so $\Fix(x)$ is a $\Lambda$-submodule of $\Dis(Q)$.

According to Proposition \ref{orb}, there is a surjection mapping $\Dis(Q)$ onto $Q_x$, given by $d \mapsto d(x)$. This surjection induces a bijection $e$ mapping $\Dis(Q)/\Fix(x)$ onto $Q_x$. To verify that $e$ is an isomorphism of quandles, notice that if $d_1,d_2 \in \Dis(Q)$ then according to Corollary \ref{dismod},
\[
e((d_1 + \Fix(x)) \triangleright (d_2 + \Fix(x))) = e(t \cdot d_1 + (1-t) \cdot d_2 + \Fix(x))
\]
\[
 = (t \cdot d_1 + (1-t) \cdot d_2)(x)=  (d_2 +t \cdot (-d_2) + t \cdot d_1)(x) 
\]
\[
= d_2 (t \cdot d_2^{-1}) (t \cdot d_1) (x)= d_2 (\beta_x d_2^{-1} \beta_x ^{-1})( \beta_x d_1 \beta_x ^{-1})(x)\]
\[
=  d_2 \beta_x d_2^{-1} d_1 \beta_x ^{-1}(x) = d_2 \beta_x d_2^{-1} d_1(x) .
\]
According to property 4 of Proposition \ref{disprop}, it follows that
\[
e((d_1 + \Fix(x)) \triangleright (d_2 + \Fix(x)))=  \beta_{ d_2(x)} d_1 (x) = d_1(x) \triangleright d_2(x) 
\]
\[
= e(d_1 + \Fix(x)) \triangleright e(d_2 + \Fix(x)).
\]
\end{proof}

A generalization of Corollary \ref{knotker} follows.

\begin{corollary}
\label{semimed}
If $Q$ is a medial quandle with only one orbit, then $Q$ is isomorphic to the standard Alexander quandle on the $\Lambda$-module $\Dis(Q)$. For any $x_0 \in Q$, a quandle isomorphism $e:\Dis(Q) \to Q$ is given by $e(d)= d(x_0)$. Moreover, if $Q$ is finitely generated, then scalar multiplication by $1-t$ defines an automorphism of $\Dis(Q)$.
\end{corollary}
\begin{proof}
Suppose $x_0 \in Q$, and $d_0 \in \Fix(x_0)$. If $x$ is any element of $Q$ then as $Q$ has only one orbit, there is a $d \in \Dis(Q)$ with $x=d(x_0)$. The group $\Dis(Q)$ is abelian, so it follows that $d_0(x)=d_0(d(x_0))=d(d_0(x_0))=d(x_0)=x$. That is, the only element of $\Fix(x_0)$ is the identity map of $Q$. Proposition \ref{orbstruc} tells us that the map $e$ is a quandle isomorphism.

It follows that the $\Lambda$-module $\Dis(Q)$ has only one orbit in its standard Alexander quandle. According to the formulas $x \triangleright y = tx + (1-t)y$ and $x \triangleright^{-1} y = t^{-1}x + (1-t^{-1}) y$, the orbit of $0$ is contained in the submodule $(1-t) \cdot \Dis(Q)$, so this submodule must be all of $\Dis(Q)$. That is, scalar multiplication by $1-t$ defines a surjective $\Lambda$-linear endomorphism $\Dis(Q) \to \Dis(Q)$.

Using the equality $\beta_z \beta_y \beta_z^{-1}=\beta_{\beta_z(y)}$ mentioned in the third paragraph of this section, it is easy to see that if $F$ is a finite generating set for the quandle $Q$, then the $\Lambda$-module $\Dis(Q)$ is generated by the elementary displacements $\beta_{f_1} \beta_{f_2}^{-1}$ with $f_1,f_2 \in F$. As $\Lambda$ is a Noetherian ring, it follows that a surjective endomorphism of $\Dis(Q)$ must be an isomorphism. \end{proof}

Another consequence of Proposition \ref{coreprop} is the following.

\begin{corollary}
\label{subq}
If $Q$ is a subquandle of a standard Alexander quandle $M$, then $Q$ is medial and semiregular, and $\Dis(Q)$ is isomorphic to the $\Lambda$-submodule of $(1-t)M$ generated by $\{(1-t)(q-q') \mid q,q' \in Q \}$.
\end{corollary}
\begin{proof}
Subquandles inherit both the medial and semiregularity properties.

Each translation $\beta_q$ of $Q$ extends to the corresponding translation $\beta_q$ of $M$. This obvious correspondence provides a homomorphism $\textup{ext}:\Dis(Q) \to \Dis(M)$, defined by $\textup{ext}(\prod \beta_{q_i}^{m_i})=\prod \beta_{q_i}^{m_i}$. Semiregularity implies that $\textup{ext}$ is well-defined, and $\textup{ext}$ is obviously injective. If $g:(1-t)M \to \Dis(M)$ is the isomorphism that appears in the proof of Proposition \ref{coreprop}, then the composition $g^{-1} \circ \textup{ext}$ maps a displacement $\prod \beta_{q_i}^{m_i} \in \Dis(Q)$ to $(1-t)(\sum m_i q_i)$. Proposition \ref{disprop} tells us that $\sum m_i = 0$; it follows that $(1-t)(\sum m_i q_i)$ is an element of the $\Lambda$-submodule of $(1-t)M$ generated by by $\{(1-t)(q-q') \mid q,q' \in Q \}$.
\end{proof}

\section{The Fundamental Medial Quandle of a Link}
\label{linkq}

Here is the definition of the fundamental medial quandle of a link, $\textup{AbQ}(L)$ in Joyce's notation \cite{J}.

\begin{definition}
\label{imq}
If $D$ is a diagram of a link $L$, then $\textup{MQ}(L)$ is the medial quandle generated by the set $\{q_a \mid a \in A(D)\}$, subject to the requirement that at each crossing $c \in C(D)$ as pictured in Fig.\ \ref{crossfig}, $q_{a_{2}} \triangleright q_{a_1}= q_{a_{3}}$.
\end{definition}

Part of Theorem \ref{main2} follows readily from Definition \ref{imq} and results of Sec.\ \ref{proof1}.

\begin{proposition}
\label{easymain2}
Suppose $L$ and $L'$ are links with $\MQ(L) \cong \MQ(L')$. Then $\Qr(L) \cong \Qr(L')$.
\end{proposition}
\begin{proof} The definition of the fundamental quandle $Q(L)$ is the same as Definition \ref{imq}, but with the word ``medial'' removed. Therefore, there is a surjective quandle map $Q(L) \to \MQ(L)$, under which the image of the $q_a$ element of $Q(L)$ is the $q_a$ element of $\MQ(L)$, for each $a \in A(D)$. 

Also, if we replace the phrase ``the medial quandle'' in Definition \ref{imq} with ``a medial quandle,'' and we replace each occurrence of $q_a$ or $q_{a_i}$ with $\gamma_D(a) \otimes 1$ or $\gamma_D(a_i) \otimes 1$, then the resulting sentence is true of $\Qr(L)$. Therefore, there is a surjective quandle map $\MQ(L) \to \Qr(L)$, under which $q_a \mapsto \gamma_D(a) \otimes 1\thickspace \allowbreak \forall a \in A(D)$. 

Theorem \ref{qpres} holds for both $Q=Q(L)$ and $Q=\Qr(L)$, so it follows that Theorem \ref{qpres} also holds for $Q=\MQ(L)$. As discussed after Theorem \ref{qpres}, we deduce that if $L$ and $L'$ are links and $f:\MQ(L) \to \MQ(L')$ is a quandle isomorphism, then after adjusting component indices in $L$ and $L'$ so that $f$ matches quandle orbits corresponding to link components with the same indices, $f$ will define an isomorphism $\Mr(L) \cong \Mr(L')$ that is compatible with the $\phi_\tau$ maps of $L$ and $L'$. That is, $L$ and $L'$ will be $\phi_\tau$-equivalent. Then Theorem \ref{main1} tells us that $\Qr(L) \cong \Qr(L')$. \end{proof}

In addition to what is stated in Proposition \ref{easymain2}, Theorem \ref{main2} asserts that the converse of Proposition \ref{easymain2} holds for knots, and fails in general. Before verifying these assertions, we present some results that will help us describe the translations and displacements of $\MQ(L)$.

\begin{proposition}
\label{imqorb}
$\textup{MQ}(L)$ has $\mu$ orbits, one for each component of $L$.
\end{proposition}
\begin{proof} The proposition follows from the corresponding properties of $Q(L)$ and $\Qr(L)$, because of the surjective quandle maps $Q(L) \to \MQ(L)$ and $\MQ(L) \to \Qr(L)$ mentioned in the proof of Proposition \ref{easymain2}. \end{proof}

A very useful notion discussed by Joyce \cite[Sec.\ 9]{J} involves describing a quandle through an ``augmentation,'' i.e.\ a group action. For $\textup{MQ}(L)$, an appropriate group is defined as follows.

\begin{definition}
\label{mg}
Let $D$ be a diagram of a link $L$. Then the \emph{medial group} $\MG(L)$ is generated by $\{g_a \mid a \in A(D)\}$, with two types of relations.
\begin{enumerate}
    \item If $a_1,a_2,a_3$ are any elements of $A(D)$ and $x_1,x_2,x_3 \in \MG(L)$ are conjugates of $g_{a_1},g_{a_2},g_{a_3}$ (respectively), then $x_1x_2^{-1}x_3=x_3x_2^{-1}x_1.$
    \item If $a_1,a_2,a_3$ are arcs appearing at a crossing of $D$ as pictured in Fig.\ \ref{crossfig}, then $g_{a_1} g_{a_2} g_{a_1}^{-1}= g_{a_3}$.
\end{enumerate}
\end{definition}

\begin{proposition}
\label{betamap}
There is a homomorphism $\beta:\MG(L) \to \Aut(\MQ(L))$ with $\beta(g_a)=\beta_{q_a} \thickspace \allowbreak \forall a \in A(D)$. 
\end{proposition}
\begin{proof}
We must show that the $\beta$ maps of $\MQ(L)$ satisfy the two kinds of relations given in Definition \ref{mg}. The relations of type 1 follow from the medial property of $\MQ(L)$; see Proposition \ref{medprop}. For the relations of type 2, suppose $a_1,a_2,a_3$ are arcs appearing at a crossing of $D$ as pictured in Fig.\ \ref{crossfig}. Then $q_{a_3}=  q_{a_2} \triangleright q_{a_1} = \beta_{q_{a_1}}(q_{a_2})$ and as noted in the third paragraph of Sec.\ \ref{struc}, it follows that $\beta_{q_{a_3}}=\beta_{q_{a_1}}\beta_{q_{a_2}}\beta^{-1}_{q_{a_1}}$.
\end{proof}

Two subsets of $\MG(L)$ will be particularly important for us.

\begin{definition}
\label{qmg}
Let $\QMG(L)$ be the subset of $\MG(L)$ that includes all conjugates of elements $g_a$, where $a \in A(D)$.
\end{definition}

\begin{definition}
\label{mg0}
Let $\MG^0(L)$ be the subset of $\MG(L)$ consisting of all products $g_{a_1}^{n_1} \cdots g_{a_k}^{n_k}$ such that $a_1, \dots, a_k \in A(D)$, $n_1, \dots, n_k \in \mathbb Z$ and $\sum n_i = 0$.
\end{definition}

The next three results show that $\QMG(L)$ is a semiregular medial quandle under conjugation. Later, we will see that $\QMG(L)$ is isomorphic to $\Qr(L)$.

\begin{proposition}
\label{qmgq}
$\QMG(L)$ is a medial quandle under conjugation: if $x,y \in \QMG(L)$, then $x \triangleright y = yxy^{-1}$. 
\end{proposition}
\begin{proof}
The first two defining properties of a quandle -- $x \triangleright x=x\thickspace \allowbreak \forall x$ and the fact that for each $y$, $\beta_y(x)=x \triangleright y$ defines a permutation -- follow from elementary properties of conjugation in groups. For the third defining property, we have
\[
(x\triangleright z) \triangleright (y \triangleright z) = zyz^{-1}zxz^{-1}(zyz^{-1})^{-1} = zy \cdot xz^{-1}zy^{-1}z^{-1}
\]
\[
 = zyxy^{-1}z^{-1}=(x\triangleright y) \triangleright z.
\]
The medial property holds in $\QMG(L)$ because 
\[
(w \triangleright x) \triangleright (y \triangleright z) = zyz^{-1}xwx^{-1}zy^{-1}z^{-1} = z(yz^{-1}x)w(x^{-1}zy^{-1})z^{-1}
\]
\[
= z(yz^{-1}x)w(yz^{-1}x)^{-1}z^{-1}= z(xz^{-1}y)w(xz^{-1}y)^{-1}z^{-1}
\]
\[
= zxz^{-1}ywy^{-1}zx^{-1}z^{-1}=(w \triangleright y) \triangleright (x \triangleright z).
\]
\end{proof}

\begin{lemma}
\label{reglem}
Suppose $n \geq 3$ is an odd integer, and $y_1,\dots,y_n \in \QMG(L)$. Then $y_1 y_2^{-1} y_3 \cdots y_{n-1}^{-1} y_n=y_n y_{n-1}^{-1} y_{n-2} \cdots y_2^{-1} y_1$ in $\MG(L)$.
\end{lemma}
\begin{proof}
If $n=3$, the assertion of the lemma is the same as part 1 of Definition \ref{mg}. Note that the $n=3$ case implies $y_1^{-1} y_2 y_3^{-1} = (y_3 y_2^{-1} y_1)^{-1}  = (y_1 y_2^{-1} y_3)^{-1} = y_3^{-1} y_2 y_1^{-1}$. If $n>3$, we have
\[
y_1 y_2^{-1} y_3 \cdots y_{n-1}^{-1} y_n=y_3 y_2^{-1} y_1 y_4^{-1} y_5\cdots y_{n-1}^{-1} y_n = \cdots
\]
\[
= y_3 y_2^{-1} y_5 y_4^{-1} \cdots y_n y_{n-1}^{-1} y_1
\]
\[
=y_3 y_4^{-1} y_5 y_2^{-1} \cdots y_n y_{n-1}^{-1} y_1 = \cdots = y_3 y_4^{-1} y_5  \cdots y_{n-1}^{-1} y_n y_2^{-1} y_1.
\]
The assertion of the lemma follows, by applying an inductive hypothesis to $y_3 y_4^{-1} y_5  \cdots y_{n-1}^{-1} y_n$. \end{proof}

\begin{proposition}
\label{imqsemi}
$\QMG(L)$ is semiregular.
\end{proposition}
\begin{proof}
Let $d \in \Dis(\QMG(L))$. Proposition \ref{disprop} implies $d=\beta_{y_1} \beta_{y_2}^{-1} \cdots \beta_{y_{2n-1}} \beta_{y_{2n}}^{-1}$ for some $y_1,\dots,y_{2n} \in \QMG(L)$. As the quandle operation of $\QMG(L)$ is given by $\beta_z(x)=zxz^{-1}$, it follows that $d$ is defined by conjugation by $y=y_1 y_2^{-1} \cdots y_{2n-1}y_{2n}^{-1}$.

If $d(x_0)=x_0$, then according to Lemma \ref{reglem}, 
\[
1=x_0^{-1} d(x_0)=x_0^{-1}yx_0y^{-1}=x_0^{-1} \cdot (y_1 y_2^{-1} \cdots y_{2n-1}y_{2n}^{-1} x_0) \cdot y_{2n} y_{2n-1}^{-1} \cdots y_2 y_1^{-1}
\]
\[
=x_0^{-1} \cdot (x_0 y_{2n}^{-1} y_{2n-1} \cdots y_2^{-1} y_1) \cdot y_{2n} y_{2n-1}^{-1} \cdots y_2 y_1^{-1}
\]
\[
= y_{2n}^{-1} y_{2n-1} \cdots y_2^{-1} y_1 \cdot y_{2n} y_{2n-1}^{-1} \cdots y_2 y_1^{-1}.
\]
Using Lemma \ref{reglem} again, we deduce that for every $x \in \QMG(L)$, 
\[
d(x)= y_1 y_2^{-1} \cdots y_{2n-1}y_{2n}^{-1} x \cdot y_{2n} y_{2n-1}^{-1} \cdots y_2 y_1^{-1}
\]
\[
= x \cdot y_{2n}^{-1} y_{2n-1} \cdots y_2^{-1} y_1 \cdot y_{2n} y_{2n-1}^{-1} \cdots y_2 y_1^{-1}=x \cdot 1 = x.
\]
That is: if $d$ has a fixed point, then $d$ is the identity map.
\end{proof}

Now we turn our attention to $\MG^0(L)$.

\begin{proposition}
\label{mgzerogen}
$\MG^0(L)$ is the subgroup of $\MG(L)$ generated by products $xy^{-1}$ such that $x,y \in \QMG(L)$.
\end{proposition}
\begin{proof}
If $H$ is the subgroup of $\MG(L)$ generated by $\{xy^{-1} \mid x,y \in \QMG(L)\}$, then certainly $H \subseteq \MG^0(L)$. We claim that $\MG^0(L) \subseteq H$.

Let $x \in \MG^0(L)$, and let $\ell(x)$ be the smallest integer $\ell$ such that there exist $x_1, \dots, x_{\ell} \in \QMG(L)$ and $\epsilon_1, \dots, \epsilon_{\ell} \in \{-1,1\}$ with $\prod x_i^{\epsilon_i}=x$ and $\sum \epsilon_i=0$. If $\ell(x)=0$, then $x=1 \in H$. The argument proceeds using induction on $\ell(x) \geq 2$. 

Case 1. If $\epsilon_1 =1 \neq \epsilon_2$, then $y=x_1x_2^{-1} \in H$, and the inductive hypothesis implies that $z=\prod_{i=3}^{\ell(x)}x_i^{\epsilon_i} \in H$. It follows that $x=yz \in H$. 

Case 2. If $\epsilon_1 =-1 \neq \epsilon_2$, then $y=x_1^{-1}x_2 = (x_1^{-1}x_2x_1)x_1^{-1} \in H$. Again, the inductive hypothesis implies that $z=\prod_{i=3}^{\ell(x)}x_i^{\epsilon_i} \in H$, and hence $x=yz \in H$.

Case 3. Suppose $\epsilon_1 = \epsilon_2$, and $j>2$ is the smallest integer with $\epsilon_j \neq \epsilon_1$. Let $g=\prod_{i=1}^{j-1}x_i^{\epsilon_i}$. Then
\[
x = \prod_{i=1}^{\ell(x)}x_i^{\epsilon_i} = gx_j^{\epsilon_j} \cdot g^{-1}  g \cdot \prod_{i=j+1}^{\ell(x)}x_i^{\epsilon_i} = (gx_jg^{-1})^{\epsilon_j} \cdot x_1^{\epsilon_1} \cdot \prod_{i=2}^{j-1}x_i^{\epsilon_i} \cdot \prod_{i=j+1}^{\ell(x)}x_i^{\epsilon_i}
\]
and since $\epsilon_1 \neq \epsilon_j$, the latter product falls under case 1 or case 2. \end{proof}

\begin{theorem}
\label{disqgen}
The following properties hold for $\MG^0(L)$.
\begin{enumerate}
    \item $\MG^0(L)$ is a normal subgroup of $\MG(L)$.
    \item The map  $\beta:\MG(L) \to \Aut(\MQ(L))$ has $\beta(\MG^0(L))=\Dis(\MQ(L))$.
    \item $\MG^0(L)$ is commutative.
    \item Let $z^*$ be a fixed element of $\QMG(L)$. Then $\MG^0(L))$ is a $\Lambda$-module, with addition given by multiplication in $\MG(L)$ and scalar multiplication given by
    \[
    t \cdot x = z^*x(z^*)^{-1}  \thickspace \allowbreak \forall x \in \MG^0(L).
    \]
    \item The module structure on $\MG^0(L)$ is independent of the choice of $z^*$. That is: if $z^*,z^{**} \in \QMG(L)$ then $z^*x(z^*)^{-1}=z^{**}x(z^{**})^{-1} \thickspace \allowbreak \forall x \in \MG^0(L)$.
    \item Let $a^* \in A(D)$ be a fixed element. Then $\MG^0(L)$ is generated, as a $\Lambda$-module, by the elements $h_a = g_ag_{a^*}^{-1}$ with $a \in A(D)$.
    \item For any crossing of $D$ as pictured in Fig. \ref{crossfig}, $h_{a_3}=(1-t)h_{a_1}+th_{a_2}$.
    \end{enumerate}
\end{theorem}
\begin{proof}
The first property follows immediately from Definition \ref{mg0}, and the second follows from the third property mentioned in Proposition \ref{disprop}.

Proposition \ref{mgzerogen} tells us that $\MG^0(L)$ is generated by elements of the form $xy^{-1}$ with $x,y \in \QMG(L)$. To verify commutativity of $\MG^0(L)$, then, it suffices to show that such elements commute:
\[
(vw^{-1})(xy^{-1})=v(w^{-1}xy^{-1})=v(yx^{-1}w)^{-1}=v(wx^{-1}y)^{-1}
\]
\[
=(vy^{-1}x)w^{-1}=(xy^{-1}v)w^{-1}=(xy^{-1})(vw^{-1}).
\]

The fourth property follows from commutativity of $\MG^0(L)$ and the fact that conjugation by $x^*$ defines an automorphism of the normal subgroup $\MG^0(L) \subset \MG(L)$. For the fifth property, note that if $z^*,z^{**} \in \QMG(L)$ then $z^*(z^{**})^{-1} \in \MG^0(L)$, so if $x \in \MG^0(L)$ then
\[
z^*x(z^{*})^{-1} = z^*(z^{**})^{-1} \cdot z^{**} x (z^{**})^{-1} \cdot z^{**}(z^*)^{-1}
\]
\[
= z^{**} x (z^{**})^{-1} \cdot z^*(z^{**})^{-1} \cdot z^{**}(z^*)^{-1}= z^{**} x (z^{**})^{-1} .
\]

To verify the sixth property, let $S$ be the $\Lambda$-submodule of $\MG^0(L)$ generated by $\{h_a \mid a \in A(D)\}$. We claim that if $a_1, \dots, a_{2n} \in A(D)$, $\epsilon_1 \dots, \epsilon_{2n} \in \{-1,1\}$ and $\sum \epsilon_i=0$, then  $x=g_{a_1}^{\epsilon_1} \dots g_{a_{2n}}^{\epsilon_{2n}} \in S$. 

Before presenting the argument, we should mention that the operation of $\MG^0(L)$ can be written using either additive or multiplicative notation. Additive notation is more natural when we think of  $\MG^0(L)$ as a $\Lambda$-module, and multiplicative notation is more natural when we think of $\MG^0(L)$ as a subgroup of $\MG(L)$. For example, we can write 
\[
g_{a_1}g_{a_2}^{-1}g_{a_3}^{-1}g_{a_4}=g_{a_1}g_{a_2}^{-1}+g_{a_3}^{-1}g_{a_4}=g_{a_1}g_{a_2}^{-1}-g_{a_4}g_{a_3}^{-1}.
\]
However, the expression $g_{a_1}g_{a_2}^{-1}g_{a_3}^{-1}+g_{a_4}$ is meaningless, because $+$ is not used for the group operation in $\MG(L)$ outside of $\MG^0(L)$. As usual, the scalar multiplication operation in a module is written using $\cdot$ or juxtaposition.

Returning to the claim, let $m$ be the number of occurrences of arcs $a_i \neq a^*$ in the list $a_1, \dots, a_{2n}$. If $m=0$, then $x=g_{a_1}^{\epsilon_1} \dots g_{a_{2n}}^{\epsilon_{2n}} = (g_{a^*})^0$ is the identity element of $\MG^0(L)$. The identity element is included in every submodule.

The argument proceeds using induction on $m>0$. If $a_1=a^*$, then 
\[
t^{-\epsilon_1} \cdot x = t^{-\epsilon_1} \cdot g_{a_1}^{\epsilon_1} \dots g_{a_{2n}}^{\epsilon_{2n}} = g_{a^*}^{-\epsilon_1}g_{a_1}^{\epsilon_1} \dots g_{a_{2n}}^{\epsilon_{2n}}g_{a^*}^{\epsilon_1}= g_{a_2}^{\epsilon_2} \dots g_{a_{2n}}^{\epsilon_{2n}}g_{a_1}^{\epsilon_1}.
\]
Repeating this as many times as necessary, we will ultimately obtain a scalar multiple $t^k \cdot x$ which is equal to a product of elements of the form $g_a^{\pm 1}$, in which the first term involves an arc other than $a^*$. As $t^k \cdot x \in S$ if and only if $x \in S$, we may as well assume that $x=g_{a_1}^{\epsilon_1} \dots g_{a_{2n}}^{\epsilon_{2n}}$ and $a_1 \neq a^*$.

If $a_1 \neq a^*$ and $\epsilon_1 =1$, then 
\[
-h_{a_1}+x=h_{a_1}^{-1}x=g_{a^*}g_{a_1}^{-1}g_{a_1}^{\epsilon_1} \dots g_{a_{2n}}^{\epsilon_{2n}} = g_{a^*}g_{a_2}^{\epsilon_2} \dots g_{a_{2n}}^{\epsilon_{2n}}.
\]
As the latter product involves only $m-1$ occurrences of arcs not equal to $a^*$, the inductive hypothesis guarantees that $-h_{a_1}+x \in S$. Of course $h_{a_1} \in S$, so it follows that $x \in S$ too.

Similarly, if $a_1 \neq a^*$ and $\epsilon_1 =-1$, then 
\[
t^{-1} \cdot h_{a_1}+x = (t^{-1} \cdot h_{a_1}) \cdot x  =g_{a^*}^{-1}g_{a_1}g_{a^*}^{-1}g_{a^*} \cdot g_{a_1}^{\epsilon_1} \dots g_{a_{2n}}^{\epsilon_{2n}}=g_{a^*}^{-1} \cdot g_{a_2}^{\epsilon_2} \dots g_{a_{2n}}^{\epsilon_{2n}} \text{,}
\]
and the latter product involves only $m-1$ occurrences of arcs not equal to $a^*$. Once again, the inductive hypothesis applies, and it tells us that $t^{-1} \cdot h_{a_1}+x \in S$. As $h_{a_1} \in S$, it follows that $x \in S$.

Turning to property 7, let $c \in C(D)$ be a crossing as pictured in Fig.\ \ref{crossfig}. Then $g_{a_3}=g_{a_1}g_{a_2}g_{a_1}^{-1}$ in $\MG(L)$, so
\[
h_{a_{3}} = g_{a_3} g_{a^*}^{-1}= g_{a_1} g_{a_2} g_{a_1}^{-1} g_{a^*}^{-1}
=  (g_{a_1} g_{a^*}^{-1}) (g_{a^*} g_{a_2} g_{a^*}^{-2}) (g_{a^*}^{2} g_{a_1}^{-1} g_{a^*}^{-1})
\]
\[
= g_{a_1} g_{a^*}^{-1}+g_{a^*} g_{a_2} g_{a^*}^{-2} + g_{a^*}^{2} g_{a_1}^{-1} g_{a^*}^{-1}= h_{a_1} +g_{a^*} (g_{a_2} g_{a^*}^{-1}) g_{a^*}^{-1} + (g_{a^*} g_{a_1} g_{a^*}^{-2})^{-1}
\]
\[
= h_{a_1} +t \cdot h_{a_2} - g_{a^*} (g_{a_1} g_{a^*}^{-1})g_{a^*}^{-1}= h_{a_1} +t \cdot h_{a_2} - t \cdot h_{a_1}.
\]
\end{proof}

Recall that if $D$ is a diagram of a link $L$, then $\fr:\Mr(L) \to \Lambda$ is the $\Lambda$-linear map with $\fr(\gamma_D(a) \otimes 1)=1\thickspace \allowbreak \forall a \in A(D)$. Of course, the kernel of $\fr$ is the submodule of $M_A^{\textup{red}}(L)$ generated by the elements $(\gamma_D(a)-\gamma_D(a')) \otimes 1$ with $a,a' \in A(D)$.

\begin{corollary}
\label{dismap}
In the situation of Theorem \ref{disqgen}, there is a $\Lambda$-linear epimorphism $e_D:\ker \fr \to \MG^0(L)$ given by $e_D((\gamma_D(a)-\gamma_D(a^*)) \otimes 1)=h_a \thickspace \allowbreak \forall a \in A(D)$.
\end{corollary}
\begin{proof}
As discussed in Sec.\ \ref{defs}, $M_A^{\textup{red}}(L)$ is the $\Lambda$-module generated by the elements $\gamma_D(a) \otimes 1$ with $a \in A(D)$, subject to the defining relations 
\[
0=\gamma_D\rho_D(c)\otimes 1 = (1-t)(\gamma_D(a_1) \otimes 1) + t(\gamma_D(a_2) \otimes 1)-(\gamma_D(a_3) \otimes 1)
\]
whenever $c \in C(D)$ is a crossing as pictured in Fig.\ \ref{crossfig}. It follows from property 7 of Theorem \ref{disqgen} that there is a $\Lambda$-linear map $M_A^{\textup{red}}(L) \to \MG^0(L)$ given by $\gamma_D(a) \otimes 1 \mapsto h_a\thickspace \allowbreak \forall a \in A(D)$. Restricting this map to $\ker \fr$ yields $e_D$. The stated formula $e_D((\gamma_D(a)-\gamma_D(a^*)) \otimes 1)=h_a$ reflects the fact that $h_{a^*}=0$ in $\MG^0(L)$.

The fact that $e_D$ is surjective follows from property 6 of Theorem \ref{disqgen}.
\end{proof}

Notice that the fixed element $a^* \in A(D)$ is helpful in stating a definition for $e_D$, and in referencing Theorem \ref{disqgen} for the fact that $e_D$ is surjective. But once we know $e_D$ is well defined, we can calculate a formula for $e_D$ that does not require $a^*$: if $a,a'$ are any elements of $A(D)$, then
\[
e_D((\gamma_D(a)-\gamma_D(a')) \otimes 1) = e_D((\gamma_D(a)-\gamma_D(a^*)) \otimes 1)-e_D((\gamma_D(a')-\gamma_D(a^*)) \otimes 1)
\]
\begin{equation}
\label{eD}
= h_a (h_{a'})^{-1} = g_a (g_{a^*})^{-1}(g_{a'} (g_{a^*})^{-1})^{-1}= g_a (g_{a'})^{-1}.
\end{equation}

We are now ready to prove the special assertion of Theorem \ref{main2} for knots.

\begin{corollary}
\label{knotcase}
If $\mu=1$, then $\MQ(L) \cong \Qr(L)$.
\end{corollary}
\begin{proof}
As noted in the proof of Proposition \ref{easymain2}, it follows from Definition \ref{imq} that for any link $L$ there is a surjective quandle map $\MQ(L) \to \Qr(L)$, with $q_a \mapsto \gamma_D(a) \otimes 1\thickspace \allowbreak \forall a \in A(D)$. As observed after Corollary \ref{dismod}, this surjective quandle map induces a surjective $\Lambda$-linear map $\Dis(\MQ(L)) \to \Dis(\Qr(L))$.

Theorem \ref{disqgen} and Corollary \ref{dismap} provide a surjective homomorphism $\beta e_D:\ker \fr \to \Dis(\MQ(L))$. The map $e_D$ is $\Lambda$-linear, and comparing item 4 of Theorem \ref{disqgen} to Corollary \ref{dismod}, we see that $\beta:\MG^0(L) \to \Dis(\MQ(L))$ is $\Lambda$-linear too. Therefore $\beta e_D:\ker \fr \to \Dis(\MQ(L))$ is surjective and $\Lambda$-linear.

Corollary \ref{subq} tells us that $\Dis(\Qr(L)) \cong (1-t) \cdot \ker \fr$, and if $\mu=1$, Corollary \ref{knotker} tells us that $(1-t) \cdot \ker \fr = \ker \fr$. Therefore, if $\mu=1$ we have surjective $\Lambda$-linear maps $\Dis(\MQ(L)) \to \Dis(\Qr(L))$ and $\Dis(\Qr(L)) \to \Dis(\MQ(L))$. As the displacement groups are finitely generated $\Lambda$-modules and $\Lambda$ is Noetherian, it follows that both surjections are isomorphisms. As each of $\Dis(\MQ(L))$ and $\Dis(\Qr(L))$ is a medial quandle with only one orbit, Corollary \ref{semimed} tells us that both $\Lambda$-module isomorphisms are also quandle isomorphisms. \end{proof}

\section{Examples for Theorem \ref{main2}}
\label{twoproof}

In this section we complete the proof of Theorem \ref{main2}, by providing a pair of links distinguished by their $\MQ$ quandles but not by their $\Qr$ quandles.

We begin with some examples of small medial quandles. Let $\id:\mathbb Z \to \mathbb Z$  be the identity map, and let $s^+,s^-:\mathbb Z \to \mathbb Z$ be the two unit shift maps, given by $s^{\pm}(n)=n \pm 1$. Let $Q$ be a set consisting of three disjoint copies of $\mathbb Z$, denoted $Q_1,Q_2$ and $Q_3$. Choose three triples $t_1,t_2,t_3$ in such a way that each triple $t_i=(t_{i1},t_{i2},t_{i3})$ has $t_{ij} \in \{\id,s^+,s^-\}\thickspace \allowbreak \forall i,j \in \{1,2,3\}$, $t_{ii}=\id\thickspace \allowbreak \forall i \in \{1,2,3\}$, and for each $j \in \{1,2,3\}$, at least one $t_{ij}$ is $s^+$ or $s^-$. Use the triples $t_i$ to define bijections $\beta_i:Q \to Q$, with $\beta_i \mid Q_j = t_{ij}:Q_j \to Q_j$. Notice that $\beta_i \beta_j = \beta_j \beta_i\thickspace \allowbreak \forall i,j \in \{1,2,3\}$. 
\begin{lemma}
\label{threeq}
Under these circumstances, $Q$ is a medial quandle under the operation $\triangleright$ defined by: if $x \in Q_j$ and $y \in Q_i$, then $x \triangleright y = \beta_i(x)$. The quandle $Q$ has three orbits, $Q_1$, $Q_2$ and $Q_3$.
\end{lemma}
\begin{proof}
If $x \in Q_i$, then $x \triangleright x= \beta_i(x) = t_{ii}(x) = \id(x) =x$. Also, each $\beta_i$ is a permutation of $Q$. If $x \in Q_k$, $y \in Q_j$ and $z \in Q_i$ then 
\[
(x\triangleright y)\triangleright z = \beta_i \beta_j(x) = \beta_j \beta_i(x) = \beta_i(x) \triangleright (y \triangleright z) =  (x \triangleright z) \triangleright (y \triangleright z)\text{,}
\]
so $Q$ is a quandle.

As $\beta_i \beta_j = \beta_j \beta_i \thickspace \allowbreak \forall i,j \in \{1,2,3\}$, the permutations $\beta_1, \beta_2, \beta_3$ generate a commutative group $G$ of permutations of $Q$. Every $\beta$ map of the quandle $Q$ is one of the three permutations $\beta_1, \beta_2, \beta_3$, so $\Dis(Q)$ is a subgroup of $G$. It follows that $\Dis(Q)$ is abelian, so $Q$ is a medial quandle. The fact that the orbits of $Q$ are $Q_1,Q_2$ and $Q_3$ follows from the hypothesis that for each $j$, at least one $t_{ij}$ is $s^+$ or $s^-$. \end{proof}

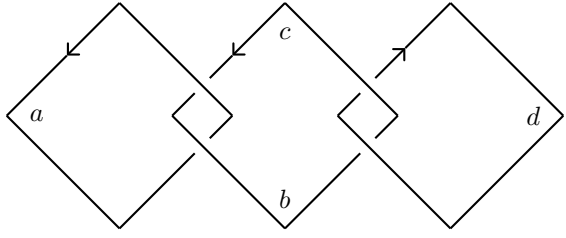
\begin{figure} [bht]
\centering
\begin{tikzpicture} [>=angle 90]
\draw [thick] [->] (0,1.5) -- (-0.7,0.8);
\draw [thick] (-0.7,0.8) -- (-1.5,0);
\draw [thick] (0,-1.5) -- (1,-0.5);
\draw [thick] (0,1.5) -- (1.5,0);
\draw [thick] (1.5,0) -- (1.2,-0.3);
\draw [thick] (-1.5,0) -- (0,-1.5);
\draw [thick] (0.7,0) -- (2.2,-1.5);
\draw [thick] (0.7,0) -- (1,0.3);
\draw [thick] [->] (2.2,1.5) -- (1.5,0.8);
\draw [thick] (1.5,0.8) -- (1.2,0.5);
\draw [thick] (2.2,1.5) -- (3.7,0);
\draw [thick] (3.7,0) -- (3.4,-0.3);
\draw [thick] (2.2,-1.5) -- (3.2,-0.5);
\draw [thick] (4.4,-1.5) -- (2.9,0);
\draw [thick] (4.4,-1.5) -- (5.9,0);
\draw [thick] (5.9,0) -- (4.4,1.5);
\draw [thick] (4.4,1.5) -- (3.8,0.9);
\draw [thick] [<-] (3.8,0.9) -- (3.4,0.5);
\draw [thick] (3.2,0.3) -- (2.9,0);

\node at (-1.1,0) {$a$};
\node at (2.2,-1.1) {$b$};
\node at (2.2,1.1) {$c$};
\node at (5.5,0) {$d$};
\end{tikzpicture}
\caption{This link is denoted $L$ in Sec.\ \ref{twoproof}.}
\label{hthreefig}
\end{figure}

Let $L$ be the link pictured in Fig.\ \ref{hthreefig}. Then $\MG(L)$ is generated by $g_a,g_b, g_c$ and $g_d$, with two types of defining relations: $xy^{-1}z=zy^{-1}x$ for all conjugates of $g_a,g_b, g_c$ and $g_d$, and four crossing relations: $g_b=g_ag_cg_a^{-1}$, $g_a=g_bg_ag_b^{-1}$, $g_d=g_cg_dg_c^{-1}$, and $g_b=g_dg_cg_d^{-1}$. Using the last relation to eliminate $g_b$, we are left with generators $g_a,g_c$ and $g_d$, and relations $g_ag_cg_a^{-1}=g_dg_cg_d^{-1}$, $g_a=g_dg_cg_d^{-1}g_ag_dg_c^{-1}g_d^{-1}$ and $g_d=g_cg_dg_c^{-1}$.

The third relation implies that $g_c$ and $g_d$ commute with each other. The first relation then implies $g_ag_cg_a^{-1}=g_c$; i.e.\ $g_a$ and $g_c$ commute with each other. It follows that $g_a$ and $g_d$ commute with each other:
\[
g_ag_d= g_a g_d \cdot g_c^{-1}g_c = g_a \cdot g_dg_c^{-1} \cdot g_c = g_a \cdot g_c^{-1} g_d \cdot g_c = g_a g_c^{-1} g_d \cdot g_c
\]
\[
= g_d g_c^{-1} g_a \cdot g_c= g_d g_c^{-1}  \cdot g_ag_c = g_d g_c^{-1}  \cdot g_cg_a = g_dg_a \text{.}
\]

Therefore, $\MG(L)$ is the free abelian group on the generators $g_a,g_c,g_d$. According to Proposition \ref{betamap}, it follows that the subgroup of $\Aut(\MQ(L))$ generated by $\{\beta_q \mid q \in \MQ(L)\}$ is commutative.

\begin{lemma}
\label{commbeta}
Let $Q$ be a quandle with $\beta_q \beta_r=\beta_r \beta_q \thickspace \allowbreak \forall q,r \in Q$. Then (a) $\beta_q = \beta_{q'}$ whenever $q$ and $q'$ are elements of the same orbit in $Q$, and (b) if $f:Q \to Q$ is a composition of $\beta$ maps of $Q$ that fixes an element $x$ of $Q$, then the restriction of $f$ to the orbit of $x$ is the identity map.
\end{lemma}
\begin{proof}
The displacement group of $Q$ is generated by products $\beta_y \beta_z^{-1}$, with $y,z \in Q$. Therefore, the hypothesis $\beta_q \beta_r=\beta_r \beta_q \thickspace \allowbreak \forall q,r \in Q$ implies that $\Dis(Q)$ is an abelian group; according to Proposition \ref{disprop}, it follows that $Q$ is a medial quandle. 

In the third paragraph of Sec.\ \ref{struc}, it was observed that the identity $\beta_{q \triangleright r} = \beta_r \beta_q \beta_r^{-1}$ holds for all elements of any quandle. In $Q$ the $\beta$ maps commute with each other, so we have $\beta_{q \triangleright r} = \beta_q \thickspace \allowbreak \forall q,r \in Q$. It follows that (a) holds.

For (b), note that if $f(x)=x$ then $f(\beta_q^n (x)) = \beta_q^n(f(x)) = \beta_q^n(x) \thickspace \allowbreak \forall q \in Q\thickspace \allowbreak \forall n \in \mathbb Z$. \end{proof}

Lemma \ref{commbeta} tells us that $\beta_{q_a},\beta_{q_c}$ and $\beta_{q_d}$ are the only $\beta$ maps of the quandle $\MQ(L)$; in particular, $\beta_{q_c} = \beta_{q_b}$. As $\beta_{q_a}(q_a)=q_a=\beta_{q_b}(q_a)$, the restrictions of $\beta_{q_a}$ and $\beta_{q_c}$ to the orbit of $\MQ(L)$ containing $q_a$ are both the identity map, so the orbit of $q_a$ is $\{\beta_{q_d}^n(q_a) \mid n \in \mathbb Z\}$. Similarly, the equalities $\beta_{q_c}(q_d)=q_d=\beta_{q_d}(q_d)$ imply that both $\beta_{q_c}$ and $\beta_{q_d}$ restrict to the identity map of the orbit of $q_d$, so this orbit is $\{\beta_{q_a}^n(q_d) \mid n \in \mathbb Z\}$. The equalities $\beta_{q_c}(q_c)=q_c= \beta_{q_d}^{-1}(q_b) = \beta_{q_d}^{-1}\beta_{q_a}(q_c)$ imply that both $\beta_{q_c}$ and $\beta_{q_d}^{-1}\beta_{q_a}$ restrict to the identity map of the orbit of $q_c$, and this orbit is $\{\beta_{q_a}^n(q_c) \mid n \in \mathbb Z\}$, with $\beta_{q_a}^n(q_c) = \beta_{q_d}^n(q_c)\thickspace \allowbreak \forall n \in \mathbb Z$.

Let $Q$ be a quandle of the type described in Lemma \ref{threeq}, corresponding to the triples $t_1 = (\id,s^+,s^+)$, $t_2=(\id,\id,\id)$ and $t_3 = (s^+,s^+, \id)$. Then the observations of the preceding paragraph imply that there is a surjective quandle map $Q \to \MQ(L)$, under which the $0$ elements of $Q_1,Q_2$ and $Q_3$ are mapped to $q_a,q_c$ and $q_d$ (respectively). (The element of $Q_2$ corresponding to $1$ is mapped to $q_b$.) As the relations required by the crossings of Fig.\ \ref{hthreefig} hold in $Q$, and $\MQ(L)$ is the largest medial quandle generated by $q_a,q_c$ and $q_d$ in which these crossing relations hold, the surjective quandle map $Q \to \MQ(L)$ must be an isomorphism.

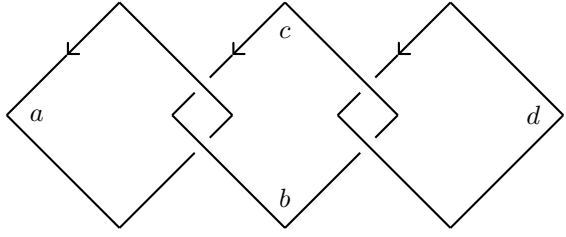
\begin{figure} [bht]
\centering
\begin{tikzpicture} [>=angle 90]
\draw [thick] [->] (0,1.5) -- (-0.7,0.8);
\draw [thick] (-0.7,0.8) -- (-1.5,0);
\draw [thick] (0,-1.5) -- (1,-0.5);
\draw [thick] (0,1.5) -- (1.5,0);
\draw [thick] (1.5,0) -- (1.2,-0.3);
\draw [thick] (-1.5,0) -- (0,-1.5);
\draw [thick] (0.7,0) -- (2.2,-1.5);
\draw [thick] (0.7,0) -- (1,0.3);
\draw [thick] [->] (2.2,1.5) -- (1.5,0.8);
\draw [thick] (1.5,0.8) -- (1.2,0.5);
\draw [thick] (2.2,1.5) -- (3.7,0);
\draw [thick] (3.7,0) -- (3.4,-0.3);
\draw [thick] (2.2,-1.5) -- (3.2,-0.5);
\draw [thick] (4.4,-1.5) -- (2.9,0);
\draw [thick] (4.4,-1.5) -- (5.9,0);
\draw [thick] (5.9,0) -- (4.4,1.5);
\draw [thick] [->] (4.4,1.5) -- (3.7,0.8);
\draw [thick] (3.7,0.8) -- (3.4,0.5);
\draw [thick] (3.2,0.3) -- (2.9,0);

\node at (-1.1,0) {$a$};
\node at (2.2,-1.1) {$b$};
\node at (2.2,1.1) {$c$};
\node at (5.5,0) {$d$};
\end{tikzpicture}
\caption{This link is denoted $L'$ in Sec.\ \ref{twoproof}.}
\label{hthreefi}
\end{figure}

Now, let $L'$ be the link pictured in Fig.\ \ref{hthreefi}. Then $\MG(L')$ is generated by $g_a,g_b, g_c$ and $g_d$, with $xy^{-1}z=zy^{-1}x$ for all conjugates of $g_a,g_b, g_c$ and $g_d$, and these four crossing relations: $g_b=g_ag_cg_a^{-1}$, $g_a=g_bg_ag_b^{-1}$, $g_d=g_cg_dg_c^{-1}$, and $g_c=g_dg_bg_d^{-1}$. Using the last relation to eliminate $g_c$, we are left with generators $g_a,g_b$ and $g_d$, and relations $g_b=g_ag_dg_bg_d^{-1}g_a^{-1}$, $g_a=g_bg_ag_b^{-1}$ and $g_d=g_dg_bg_d^{-1}g_dg_dg_b^{-1}g_d^{-1}$.

The second relation implies that $g_a$ and $g_b$ commute with each other. The first implies $g_a^{-1}g_bg_a=g_dg_bg_d^{-1}$; as $g_a$ and $g_b$ commute with each other, this implies $g_b=g_dg_bg_d^{-1}$, so $g_b$ and $g_d$ commute with each other. Then 
\[
g_ag_d= g_a g_d \cdot g_b^{-1}g_b = g_a \cdot g_dg_b^{-1} \cdot g_b = g_a \cdot g_b^{-1} g_d \cdot g_b = g_a g_b^{-1} g_d \cdot g_b
\]
\[
= g_d g_b^{-1} g_a \cdot g_b= g_d g_b^{-1}  \cdot g_ag_b = g_d g_b^{-1}  \cdot g_bg_a = g_dg_a \text{,}
\]
so $g_a$ and $g_d$ commute with each other. That is, $\MG(L')$ is commutative. It follows that $\MG(L')$ is the free abelian group on the generators $g_a,g_b,g_d$. According to Proposition \ref{betamap}, all of the $\beta$ maps of $\MG(L')$ commute with each other, and according to Lemma \ref{commbeta}, the only $\beta$ maps of the quandle $\MQ(L')$ are $\beta_{q_a},\beta_{q_b}=\beta_{q_c}$ and $\beta_{q_d}$. The equalities $\beta_{q_a}(q_a)=q_a=\beta_{q_b}(q_a)$ and $\beta_{q_c}(q_d)=q_d=\beta_{q_d}(q_d)$ imply that the restrictions of $\beta_{q_a}$ and $\beta_{q_c}$ to the $q_a$ orbit of $\MQ(L')$ are both the identity map, and the restrictions of $\beta_{q_c}$ and $\beta_{q_d}$ to the $q_d$ orbit are both the identity map. The equalities $\beta_{q_c}(q_c)=q_c= \beta_{q_d}(q_b) = \beta_{q_d}\beta_{q_a}(q_c)$ imply that both $\beta_{q_c}$ and $\beta_{q_d}\beta_{q_a}$ restrict to the identity map of the $q_c$ orbit.

Let $Q'$ be the quandle from Lemma \ref{threeq} corresponding to the triples $t'_1=(\id,s^+,s^+)$, $t'_2=(\id,\id,\id)$ and $t'_3=(s^+,s^-,\id)$. Then the discussion above implies that there is a surjective quandle map $Q' \to \MQ(L')$, under which the $0$ elements of $Q'_1,Q'_2,Q'_3$ are mapped to $q_a,q_c$ and $q_d$ respectively. The crossing relations from Fig.\ \ref{hthreefi} hold in $Q'$, so this quandle map is an isomorphism. 

Notice that the quandle $\MQ(L)$ has these properties: there are only three $\beta$ maps, one from each orbit; one $\beta$ map is the identity; and the other two $\beta$ maps have the same restriction to the orbit whose $\beta$ map is the identity. In contrast, $\MQ(L')$ has these properties: there are only three $\beta$ maps, one from each orbit; one $\beta$ map is the identity; and the restrictions of the other two $\beta$ maps to the orbit whose $\beta$ map is the identity are inverses of each other. It is clear that $\MQ(L)$ and $\MQ(L')$ are not isomorphic.

To complete the proof of Theorem \ref{main2}, we need to show that $\Qr(L)$ and $\Qr(L')$ are isomorphic. According to Theorem \ref{main1}, it is enough to show that $L$ and $L'$ are $\phi_\tau$-equivalent. In fact, $L$ and $L'$ are related through the stronger multivariate version of  $\phi_\tau$-equivalence, which we call \emph{Crowell equivalence}.

\begin{proposition}
There is a $\Lambda_\mu$-linear isomorphism $f:M_A(L) \to M_A(L')$, with $\phi_L = \phi_{L'} \circ f:M_A(L) \to I_\mu$.
\end{proposition}

\begin{proof}
We refer to Sec.\ \ref{defs} for definitions.

Abusing notation, we use $D$ for the link diagrams in both Fig.\ \ref{hthreefig} and Fig.\ \ref{hthreefi}, and we use $a,b,c,d$ for the arcs in both diagrams. The components of $L$ and $L'$ are indexed in order, from left to right in Figs.\ \ref{hthreefig} and \ref{hthreefi}.

The two crossings on the left in Figs.\ \ref{hthreefig} and \ref{hthreefi} provide the relations $\gamma_D(b) =  (1- t_2)\gamma_D(a) + t_1\gamma_D(c)$ and $\gamma_D(a)= (1-t_1)\gamma_D(b) + t_2\gamma_D(a)$, in both $M_A(L)$ and $M_A(L')$. The second relation is equivalent to $(1-t_2)\gamma_D(a)= (1-t_1)\gamma_D(b)$, and with this equality, the first relation is equivalent to $\gamma_D(b)=\gamma_D(c)$.

Keeping in mind that $\gamma_D(b)=\gamma_D(c)$ in both $M_A(L)$ and $M_A(L')$, the two crossings on the right in Figs.\ \ref{hthreefig} and \ref{hthreefi} provide the same two relations in $M_A(L)$ and $M_A(L')$: $\gamma_D(d)=(1-t_3)\gamma_D(b) + t_2 \gamma_D(d)$ and $\gamma_D(b)=(1-t_2)\gamma_D(d) + t_3 \gamma_D(b)$. Both of these relations are equivalent to $(1-t_3)\gamma_D(b) = (1-t_2) \gamma_D(d)$.

In summary, $M_A(L)$ and $M_A(L')$ are both generated by $\gamma_D(a)$, $\gamma_D(b)$ and $\gamma_D(d)$, subject to the relations $(1-t_2) \gamma_D(a)=(1-t_1)\gamma_D(b) $ and $(1-t_3)\gamma_D(b) = (1-t_2) \gamma_D(d)$. The obvious isomorphism $f:M_A(L) \to M_A(L')$ satisfies the statement.
\end{proof}

To describe the $\Lambda$-modules $\Mr(L)$ and $\Mr(L')$, we modify the descriptions of $M_A(L)$ and $M_A(L')$ by replacing each indeterminate $t_i$ with $t$, and replacing each element $\gamma_D(x)$ with $\gamma_D(x) \otimes 1$. We conclude that
\[
\Mr(L) \cong \Mr(L') \cong \Lambda \oplus (\Lambda/(1-t)) \oplus (\Lambda/(1-t)) \text{,}
\]
with the direct summands generated by $\gamma_D(b) \otimes 1$, $(\gamma_D(a)-\gamma_D(b)) \otimes 1$ and $(\gamma_D(d)-\gamma_D(b)) \otimes 1$, respectively. It is apparent that $L$ and $L'$ are $\phi_\tau$-equivalent.

\section{Longitudes}
\label{longs}

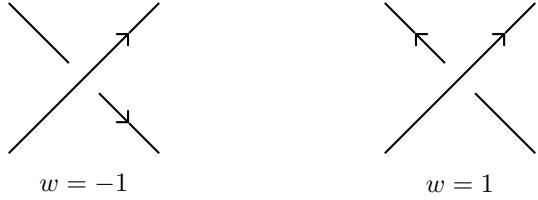
\begin{figure} [bth]
\centering
\begin{tikzpicture} [>=angle 90]
\draw [thick] (1,1) -- (0.6,0.6);
\draw [thick] [<-] (0.6,0.6) -- (-1,-1);
\draw [thick] (-1,1) -- (-.2,0.2);
\draw [thick] [->] (0.2,-0.2) -- (0.6,-0.6);
\draw [thick] (0.6,-0.6) -- (1,-1);
\draw [thick] (6,1) -- (5.6,0.6);
\draw [thick] [<-] (5.6,0.6) -- (4,-1);
\draw [thick] (4,1) -- (4.4,0.6);
\draw [thick] [<-] (4.4,0.6) -- (4.8,0.2);
\draw [thick] (5.2,-0.2) -- (6,-1);
\node at (0,-1.4) {$w=-1$};
\node at (5,-1.4) {$w=1$};
\end{tikzpicture}
\caption{The writhe of a crossing is denoted $w$.}
\label{writhefig}
\end{figure}

We would like to have more information about the map $e_D$. To obtain this information it will be useful to consider a special type of link diagram, defined using the familiar notion of writhe. See Fig. \ref{writhefig}.

\begin{definition}
\label{wbal}
Let $D$ be a link diagram. Then $D$ has \emph{alternating writhes} if every arc $a \in A(D)$ occurs as the underpassing arc of two crossings, one of writhe $-1$ and the other of writhe $1$.
\end{definition}

\begin{proposition}
\label{altw}
Every classical link has a diagram with alternating writhes.
\end{proposition}
\begin{proof}
Start with any diagram $D$ of $L$. If $L$ has a component $K_i$ which is not the underpassing component of any crossing of $D$, insert a trivial crossing into the one arc of $D$ that represents $K_i$. (Trivial crossings are pictured in Fig.\ \ref{trivf}.) We now have a diagram $D'$ in which every arc appears as the underpassing arc of at least one crossing.

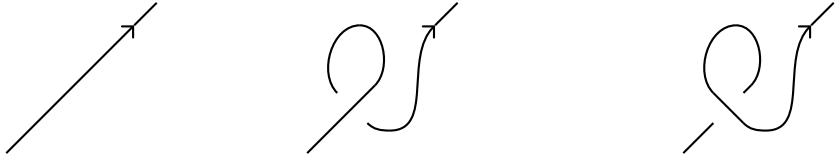
\begin{figure} [bth]
\centering
\begin{tikzpicture} [>=angle 90]
\draw [thick] [->] (-5,-1) -- (-3.3,0.7);
\draw [thick](-3,1) -- (-3.3,0.7);
\draw [thick] (1,1) -- (0.7,0.7);
\draw [thick] (-0.1,-0.1) -- (-1,-1);
\draw [thick] (-0.2,-0.6) to [out=-45, in=180] (0.1,-0.7);
\draw [thick] [->] (0.1,-0.7) to [out=0, in=-135] (0.7,0.7);
\draw [thick] (-0.1,-0.1) to [out=45, in=0] (-0.3, 0.7);
\draw [thick] (-0.3, 0.7) to [out=180, in=135] (-0.6, -0.2);
\draw [thick] (6,1) -- (5.7,0.7);
\draw [thick] (5-0.6,-0.6) -- (5-1,-1);
\draw [thick] (5-0.2,-0.6) to [out=-45, in=180] (5.1,-0.7);
\draw [thick] [->] (5.1,-0.7) to [out=0, in=-135] (5.7,0.7);
\draw [thick] (5-0.1,-0.1) to [out=45, in=0] (5-0.3, 0.7);
\draw [thick] (5-0.3, 0.7) to [out=180, in=135] (5-0.6, -0.2);
\draw [thick] (5-0.6, -0.2) -- (5-0.2,-0.6);
\draw [thick] (5-0.2,-0.2) -- (5-0.1, -0.1);
\end{tikzpicture}
\caption{A trivial crossing of either writhe may be inserted into an arc.}
\label{trivf}
\end{figure}

For every arc $a \in A(D')$ that appears as the underpassing arc only at crossing(s) of one writhe value $w$, insert a trivial crossing of writhe $-w$ into $a$. The effect is to split $a$ into two distinct arcs, each of which appears as the underpassing arc at two crossings of opposite writhe.
\end{proof}

Now, suppose $L$ is a link, and $D$ is a diagram of $L$ with alternating writhes. For each component $K_i$ of $L$, we choose an arbitrary arc $b_{i0} \in A(D)$ with $\kappa_D(b_{i0})=i$, and we start walking along $b_{i0}$, in the direction given by the orientation of $K_i$. When we reach the end of $b_{i0}$, we index that crossing as $c_{i0}$, the overpassing arc at that crossing as $a_{i0}$, and the next arc of $K_i$ as $b_{i1}$. The alternating writhes condition guarantees that as we walk along $K_i$, we pass under an even number of crossings. By the time we get back to $b_{i0}$, we will have indexed these crossings as $c_{i0}, \dots, c_{i(2n_i-1)}$, indexed the arcs of $K_i$ as $b_{i0}, \dots, b_{i(2n_i-1)}$, and indexed the overpassing arcs of the crossings over $K_i$ as $a_{i0}, \dots, a_{i(2n_i-1)}$. We consider the second indices of $a_{ij},b_{ij}$ and $c_{ij}$ modulo $2n_i$.

\newpage

Let $j \in \{1, \dots, 2n_i \}$, and let $w_{ij}$ be the writhe of the crossing $c_{ij}$. We have the following equalities in $\Mr(L)$ and $\MQ(L)$ (respectively).
\begin{align}
\label{congruence}
\gamma_D(b_{ij}) \otimes 1 = (1-t^{-w_{ij}})(\gamma_D(a_{ij}) \otimes 1)+t^{-w_{ij}}(\gamma_D(b_{i(j+1)}) \otimes 1)
\end{align}
\begin{align}
\label{congruencea}
q_{b_{ij}} = q_{b_{i(j+1)}} \triangleright^{-w_{ij}} q_{a_{ij}} =  \beta_{q_{a_{ij}}}^{-w_{ij}}(q_{b_{i(j+1)}})
\end{align}

\begin{definition}
\label{longitudes}
Under these circumstances, for each $i \in \{1, \dots, \mu \}$ the $i^{th}$ \emph{longitude} of $L$ is
\[
\chi_i=\sum_{j=1}^{2n_i} (-w_{ij})(\gamma_D(a_{ij}) \otimes 1)  \in \ker \fr .
\]
\end{definition}

\begin{proposition}
\label{trivdis}
For each $i \in \{1, \dots, \mu \}$, $\chi_i$ has the following properties.
\begin{enumerate}
    \item  $(1-t) \cdot \chi_i=0$.
    \item The image of $\chi_i$ under the map $e_D:\ker \fr \to \MG^0(L)$ is an element of the center of $\MG(L)$.
    \item The image of $\chi_i$ under the composition $\beta e_D:\ker \fr \to \Aut(\MQ(L))$ is a displacement of $\MQ(L)$ whose restriction to the $K_i$ orbit of $\MQ(L)$ is the identity map.
    \item If $\mu=1$, then $\chi_1 = 0$.
\end{enumerate}
\end{proposition}
\begin{proof} For convenience, we assume that the arcs and crossings of $D$ have been indexed so that $w_{i1}=-1$. $D$ has alternating writhes, so $w_{ij}$ is always $(-1)^j$.

According to (\ref{congruence}), if $i \in \{1, \dots, \mu \}$ and $j \in \{1, \dots, 2n_i\}$ then we have
\[
(\gamma_D(b_{ij}) \otimes 1)  = (1-t^{-w_{ij}})(\gamma_D(a_{ij}) \otimes 1) + t^{-w_{ij}}(\gamma_D(b_{i(j+1)}) \otimes 1).
\]
If $j$ is odd, it follows that 
\[
(1-t) \cdot (-w_{ij}(\gamma_D(a_{ij}) \otimes 1))=(\gamma_D(b_{ij}) \otimes 1) -t(\gamma_D(b_{i(j+1)}) \otimes 1).
\]
If $j$ is even, it follows that 
\[
(1-t^{-1})(\gamma_D(a_{ij}) \otimes 1)=(\gamma_D(b_{ij}) \otimes 1) -t^{-1}(\gamma_D(b_{i(j+1)}) \otimes 1)
\]
and hence
\[
(1-t) \cdot (-w_{ij}(\gamma_D(a_{ij}) \otimes 1))= t(\gamma_D(b_{ij}) \otimes 1)-(\gamma_D(b_{i(j+1)}) \otimes 1).
\]
Therefore,
\begin{equation*}
\begin{split}
(1-t) \cdot \chi_i = & \sum_{j=1}^{2n_i}(1-t) \cdot (-w_{ij}(\gamma_D(a_{ij}) \otimes 1))
\\
= & \sum_{\text{odd }j} ((\gamma_D(b_{ij}) \otimes 1) -t(\gamma_D(b_{i(j+1)}) \otimes 1)) 
\\
& + \sum_{\text{even } j} (t(\gamma_D(b_{ij}) \otimes 1)-(\gamma_D(b_{i(j+1)}) \otimes 1)).
\end{split}
\end{equation*}
The total is $0$ because when $j$ is odd, $\gamma_D(b_{ij}) \otimes 1$ occurs with coefficient $1$ in the penultimate sum and coefficient $-1$ in the last sum, and when $j$ is even, $\gamma_D(b_{ij}) \otimes 1$ occurs with coefficient $-t$ in the penultimate sum and coefficient $t$ in the last sum.

For the second property, recall that $e_D$ is $\Lambda$-linear by Corollary \ref{dismap}. Therefore 
\[
e_D(\chi_i) = (1-t+t) \cdot e_D( \chi_i)=e_D((1-t) \cdot \chi_i)+t \cdot e_D( \chi_i)=0+t \cdot e_D(\chi_i) = t \cdot e_D(\chi_i).
\]
According to Theorem \ref{disqgen}, it follows that for every $z^* \in \QMG(L)$,
\[
e_D(\chi_i) = t \cdot e_D(\chi_i) = z^* \cdot e_D(\chi_i) \cdot (z^*)^{-1}.
\]
That is, $e_D(\chi_i)$ commutes with every  $z^* \in \QMG(L)$. The elements of $\QMG(L)$ generate the group $\MG(L)$, so $e_D(\chi_i)$ commutes with every element of $\MG(L)$.

For the third property, notice that formula (\ref{eD}) implies
\[
e_D(\chi_i)= \sum_{k=1}^{n_i} e_D((\gamma_D(a_{i(2k-1)}) \otimes 1)-(\gamma_D(a_{i(2k)}) \otimes 1) ) 
\]
\[
= \prod_{k=1}^{n_i} (g_{a_{i(2k-1)}}g_{a_{i(2k)}}^{-1})=g_{a_{i1}}^{-w_{i1}}g_{a_{i2}}^{-w_{i2}} \cdots g_{a_{i(2n_i)}}^{-w_{i(2n_i)}}.
\]
According to Proposition \ref{betamap}, $\beta:\MG(L) \to \Aut(\MQ(L))$ is the homomorphism with $\beta(g_a)=\beta_{q_a}\thickspace \allowbreak \forall a \in A(D)$. As $b_{i(2n_i+1)}=b_{i1}$, the relations (\ref{congruencea}) imply that 
\[
\beta e_D(\chi_i) (q_{b_{i1}}) = \beta(g_{a_{i1}})^{-w_{i1}}\beta(g_{a_{i2}})^{-w_{i2}} \cdots \beta(g_{a_{i(2n_i)}})^{-w_{i(2n_i)}}(q_{b_{i(2n_i+1)}})
\]
\[
=\beta(g_{a_{i1}})^{-w_{i1}}\beta(g_{a_{i2}})^{-w_{i2}} \cdots \beta(g_{a_{i(2n_i-1)}})^{-w_{i(2n_i-1)}}(q_{b_{i(2n_i)}})
\]
\[
=  \cdots = \beta(g_{a_{i1}})^{-w_{i1}}(q_{b{i2}}) = q_{b{i1}} .
\]

The group $\Dis(\MQ(L))$ is commutative, so every $d \in \Dis(\MQ(L))$ has $\beta e_D(\chi_i) (d(q_{b_{i1}}))=d(\beta e_D(\chi_i) (q_{b_{i1}}))=d(q_{b_{i1}})$. According to Proposition \ref{orb}, every element of the $K_i$ orbit of $\MQ(L)$ is $d(q_{b{i1}})$ for some $d \in \Dis(\MQ(L))$. It follows that $\beta e_D(\chi_i)$ fixes every element of the $K_i$ orbit of $\MQ(L)$. 

The fourth property follows immediately from the first property and Corollary \ref{knotker}. \end{proof}

\begin{proposition}
\label{fed}
Let $D$ be a diagram of $L$, with alternating writhes. Then $\ann(1-t) = \{ x \in M_A^{\textup{red}}(L) \mid (1-t)x=0\}$ is the submodule of $\Mr(L)$ generated by $\chi_1, \dots, \chi_\mu$.
\end{proposition}
\begin{proof}
The first property of Proposition \ref{trivdis} implies that $\chi_1, \dots, \chi_\mu \in \ann(1-t)$. 

Now, suppose that $x \in \Mr(L)$ has $(1-t)x=0$. We must show that $x$ is equal to a $\Lambda$-linear combination of $\chi_1, \dots, \chi_\mu$.
Recall the exact sequence 
\begin{equation*}
\Lambda_{\mu}^{C(D)} \xrightarrow{\rho_D} \Lambda_{\mu}^{A(D)} \xrightarrow{\gamma_D} M_A(L) \to 0 \textup{,}
\end{equation*}
discussed in Sec.\ \ref{defs}. If $\textup{id}$ denotes the identity map of $\Lambda$, then the right exactness of tensor products yields an exact sequence
\begin{equation*}
\Lambda_{\mu}^{C(D)} \otimes_{\Lambda_{\mu}} \Lambda \xrightarrow{\rho_D \otimes \textup{id}} \Lambda_{\mu}^{A(D)}  \otimes_{\Lambda_{\mu}} \Lambda \xrightarrow{\gamma_D \otimes \textup{id}} M_A^{\textup{red}}(L) \to 0 \textup{.}
\end{equation*}
Therefore, there is a function $f_x:A(D) \to \Lambda$ such that $x= (\gamma_D \otimes \textup{id})(x')$, where
\begin{equation}
\label{lastprop1}
 x'=\sum_{a \in A(D)} f_x(a) (a \otimes 1).
\end{equation}
As $(1-t)x' \in \ker (\gamma_D \otimes \textup{id})$, there is also a function $g_x:C(D) \to \Lambda$ such that
\begin{equation}
\label{lastprop2}
(1-t)x' = \sum_{c \in C(D)} g_x(c) \cdot (\rho_D(c) \otimes 1) \textup{.}
\end{equation}

The functions $f_x$ and $g_x$ are not unique. If $c \in C(D)$ then for any element $\lambda_c \in \Lambda$, we may add $\lambda_c \cdot (\rho_D(c) \otimes 1)$ to the sum on the right-hand side of (\ref{lastprop1}) without changing the fact that $x = (\gamma_D \otimes \textup{id})(x')$, and (\ref{lastprop2}) will remain valid so long as $(1-t)\lambda_c$ is added to $g_x(c)$. In particular, for each $c \in C(D)$ there is a $\lambda_c \in \Lambda$ such that $(1-t) \lambda_c = \epsilon g_x(c) - g_x(c)$. The result of adding $\epsilon g_x(c) - g_x(c)$ to $g_x(c)$ is to replace $g_x(c)$ with $\epsilon g_x(c) - g_x(c) + g_x(c) = \epsilon g_x(c)$, which is an integer. It follows that we may assume $g_x(c) \in \mathbb Z\thickspace \allowbreak \forall c \in C(D)$, without loss of generality.

We now claim that for every arc $a \in A(D)$, the values of $g_x(c)$ for the two crossings at which $a$ is an underpassing arc are negatives of each other. To see why the claim is true, notice first that (\ref{lastprop1}) and (\ref{lastprop2}) yield
\begin{equation}
\label{lastprop3}
\sum_{a \in A(D)} (1-t)f_x(a) (a \otimes 1) = \sum_{c \in C(D)} g_x(c) \cdot (\rho_D(c) \otimes 1) \textup{.}
\end{equation}
This equality holds in the $\Lambda$-module $\Lambda_{\mu}^{A(D)}  \otimes_{\Lambda_{\mu}} \Lambda$, which is freely generated by the elements $a \otimes 1$ with $a \in A(D)$. Therefore, for each $a \in A(D)$ the coefficients of $a \otimes 1$ on the two sides of (\ref{lastprop3}) are precisely equal. 

Suppose $c \in C(D)$ has $g_x(c) \neq 0$, and $a \in A(D)$ is the arc that corresponds to $a_2$, when $c$ is pictured as in Fig.\ \ref{crossfig}. Then the contribution of the term $g_x(c) \cdot (\rho_D(c) \otimes 1)$ to the coefficient of $a \otimes 1$ on the right-hand side of (\ref{lastprop3}) is $g_x(c) \cdot t$. Let $c'$ be the other crossing of $D$ at which $a$ is one of the underpassing arcs. It is easy to see that the alternating writhes property guarantees that $a$ plays the same role at $c'$, i.e., $a_2$ rather than $a_3$, as pictured in Fig.\ \ref{crossfig}. Therefore the contribution of the term $g_x(c') \cdot (\rho_D(c') \otimes 1)$ to the coefficient of $a \otimes 1$ on the right-hand side of (\ref{lastprop3}) is $g_x(c') \cdot t$. Aside from $c$ and $c'$, this arc $a$ is incident only at crossings $c''$ where it is the overpassing arc, and for such a crossing $c''$, the contribution of the term $g_x(c'') \cdot (\rho_D(c'') \otimes 1)$ to the coefficient of $a \otimes 1$ on the right-hand side of (\ref{lastprop3}) is divisible by $1-t$. The coefficient of $a \otimes 1$ on the left-hand side of (\ref{lastprop3}) is divisible by $1-t$, so it follows that $g_x(c)=-g_x(c')$, as claimed.

The argument for an arc $a$ that plays the role of $a_3$ in Fig.\ \ref{crossfig} is almost the same. The only difference is that the contributions from $c$ and $c'$ are $g_x(c) \cdot (-1)$ and $g_x(c') \cdot (-1)$, rather than $g_x(c) \cdot t$ and $g_x(c') \cdot t$. This completes the proof of the claim.

As $D$ has alternating writhes, the claim can also be stated as follows. For each arc $a \in A(D)$, there is an integer $m_a$ such that the two crossings at which $a$ is an underpassing arc both satisfy the equality $g_x(c) = w(c)m_a$. At each crossing there is only one value of $g_x(c)$, so the two underpassing arcs must have the same value of $m_a$. Walking from crossing to crossing along the arcs of $D$, we deduce that the value of $m_a$ is constant on each component $K_i$ of $L$. We denote this constant value $m_i$. 

While proving the claim we showed that on the right-hand side of (\ref{lastprop3}), all of the contributions from underpassing arcs cancel each other. This leaves only the contributions from overpassing arcs. That is, if we recall the indexing convention for $a_{ij},b_{ij},c_{ij}$ mentioned after Proposition \ref{altw}, then 
\[
\sum_{a \in A(D)} (1-t)f_x(a) (a \otimes 1) =\sum_{c \in C(D)} g_x(c) \cdot (\rho_D(c) \otimes 1) 
\]
\[
= \sum_{i=1}^\mu \sum_{j=1}^{2n_i}g_x(c_{ij}) \cdot (1-t) (a_{ij} \otimes 1)= \sum_{i=1}^\mu \sum_{j=1}^{2n_i}w(c_{ij})m_i \cdot (1-t) (a_{ij} \otimes 1).
\]

Once again, this equality holds in the free $\Lambda$-module $\Lambda_\mu ^{A(D)} \otimes_{\Lambda_\mu} \Lambda$, so for every $a \in A(D)$, the coefficients of $a \otimes 1$ in the first and last displayed sums must be precisely equal. It follows that the factors of $1-t$ may be canceled, so
\[
x' = \sum_{a \in A(D)} f_x(a) (a \otimes 1)= \sum_{i=1}^\mu \sum_{j=1}^{2n_i}w(c_{ij})m_i \cdot (a_{ij} \otimes 1)
\]
\[
= -\sum_{i=1}^\mu m_i \sum_{j=1}^{2n_i}(-w(c_{ij})) \cdot (a_{ij} \otimes 1)
\]
and hence $x= (\gamma_D \otimes \textup{id})(x') = - \sum m_i \chi_i$.
\end{proof}

\section{$\Qr(L)$ and $\QMG(L)$}
\label{threeproof}

Our final result is that the quandles $\Qr(L)$ and $\QMG(L)$ are always isomorphic. The first part of the proof involves the Alexander module of the group $\MG(L)$. We give a brief summary of the theory regarding this module, and refer to Crowell \cite{C3} for a thorough account.

The integral group ring $\mathbb Z (\MG(L))$ consists of formal linear combinations $\sum n_i g_i$, where the $n_i$ are integers and the $g_i$ are elements of $\MG(L)$. The integral group ring $\mathbb Z (\MG(L)/\MG(L)')$ of the abelianization of $\MG(L)$ is defined analogously, and it is made into a $\mathbb Z (\MG(L))$-module using the multiplication of $\mathbb Z (\MG(L))$. The \emph{augmentation ideal} $I \subset \mathbb Z (\MG(L))$ is the ideal generated by the elements $g-1$, where $g \in \MG(L)$. The \emph{Alexander module} of $\MG(L)$ is the tensor product 
\[
\mathbb Z (\MG(L) / \MG(L)') \otimes _{\mathbb Z (\MG(L))} I = M \text{,}
\]
considered as a $\mathbb Z (\MG(L)/\MG(L)')$-module via multiplication in the first factor of the tensor product. That is, if $\alpha:\MG(L) \to \MG(L)/\MG(L)'$ is the canonical map onto the quotient, $g,h \in \MG(L)$ and $i \in I$, then 
\[
\alpha(g) \cdot (\alpha(h) \otimes i)=  \alpha(gh) \otimes i= \alpha(g) \otimes (hi)  = 1 \otimes (ghi).
\]

Definition \ref{mg} implies that  $\MG(L)/\MG(L)'$ is a free abelian group of rank $\mu$, with one generator for each component $K_i$ of $L$. The generator corresponding to $K_i$ is $\alpha(g_a)$, for all $a \in A(D)$ with $\kappa_D(a)=i$. There is then a natural isomorphism between $\mathbb Z (\MG(L) / \MG(L)')$ and $\Lambda_\mu$, under which the generator of $\mathbb Z (\MG(L) / \MG(L)')$ corresponding to $K_i$ is mapped to $t_i$.

As explained by Crowell \cite{C3}, the Alexander module $M$ is a finitely presented $\Lambda_\mu$-module. The generators in the presentation are the elements $1 \otimes (g_a - 1)$, where $a \in A(D)$. The relations in the presentation of the module $M$ are obtained by taking the free derivatives of the relators in the presentation of the group $\MG(L)$ given in Definition \ref{mg}, and then applying the abelianization map $\alpha$. 

We recall the definition of the free derivatives. Suppose $F$ is the free group on the set $X$, and 
\[
w = \prod\limits _{i=1}^n x_i^{\epsilon_i} \in F \text{,}
\]
where $x_1, \dots, x_n \in X$ and $\epsilon_1, \dots , \epsilon_n \in \{-1, 1\}$. For $1 \leq i \leq n$, define the $i$th initial segment of $w$ as follows:
\[
w_i = \begin{cases}
\prod\limits_{j=1}^{i-1} x_j^{\epsilon_j} , & \text{if } \epsilon_i=1 \\
\prod\limits _{j=1}^{i} x_j^{\epsilon_j} , & \text{if } \epsilon_i=-1.
\end{cases}
\]
Then for each $x \in X$, the \emph{free derivative} of $w$ with respect to $x$ is
\[
\label{fdc}
\frac{\partial w}{\partial x} = \sum\limits_{x_i=x} \epsilon_i w_i \in \mathbb Z F.
\]

Recall that there are two types of relators in Definition \ref{mg}: if $b,c,d$ are conjugates of $g_a$ elements then there is a relator $bc^{-1}db^{-1}cd^{-1}$, and if $a_1,a_2,a_3$ are arcs that appear at a crossing of $D$ as pictured in Fig.\ \ref{crossfig}, then there is a relator $g_{a_1} g_{a_2} g_{a_1}^{-1} g_{a_3}^{-1}$.
\begin{theorem}
\label{mgmod}
If $M$ is the Alexander module of the group $\MG(L)$, then there is an isomorphism of $\Lambda$-modules
\[
\widehat \gamma_D:M \otimes _{\Lambda_\mu}\Lambda \to \Mr(L) \text{,}
\]
with $\widehat \gamma_D((1 \otimes (g_a-1))\otimes 1)=\gamma_D(a) \otimes 1 \thickspace \allowbreak \forall a \in A(D)$.
\end{theorem}
\begin{proof}
Right exactness of tensor products implies that a presentation of the $\Lambda$-module $M \otimes _{\Lambda_\mu}\Lambda$ can be obtained from the presentation of the $\Lambda_\mu$-module $M$ described above, by applying $\tau$ to all coefficients.  

We claim that for a relator of the form $r=bc^{-1}db^{-1}cd^{-1}$, all of the resulting $\tau$ values are $0$. There are several different places in the relator where a generator might appear; we consider two of them, and leave it to the reader to consider the rest. Suppose first that $c$ is a conjugate of $g_a$, say $c=eg_ae^{-1}$. This appearance of $g_a$ in the middle of $c$ contributes two terms to the free derivative $\partial r / \partial g_a$, namely, $-beg_a^{-1}$ and $bc^{-1}db^{-1}e$. Under the composition of the abelianization map $\alpha:\MG(L) \to \MG(L)/\MG(L)'$ and the map $\tau: \mathbb Z(\MG(L)/\MG(L)') \cong \Lambda_\mu \to \Lambda$, each of $b,c,d, g_a$ is mapped to $t$. Hence the image of $-beg_a^{-1}+bc^{-1}db^{-1}e$ is the same as the image of $-e+e$; of course, this image is $0$. For another example, suppose there is an appearance of $g_a$ in $c$ that is not in the middle of $c$; say $c=e_1g_a^{-1}e_2g_{a'}e_2^{-1}g_ae_1^{-1}$. These two appearances of $g_a$ in $c$ also provide two corresponding appearances of $g_a$ in $c^{-1}$, and the total contribution to $\partial r/ \partial g_a$ of these four appearances of $g_a$ is
\[
-be_1g_a^{-1}  + be_1g_a^{-1}e_2g_{a'}^{-1}e_2^{-1} - bc^{-1}db^{-1}e_1g_a^{-1}+bc^{-1}db^{-1}e_1g_a^{-1}e_2 g_{a'}e_2^{-1}.
\]
As $b,c,d,g_a $ and $g_{a'}$ are all mapped to $t$, the image of this total in $\Lambda$ is the same as the image of 
\[
-e_1 + e_1 g_a^{-1} -e_1g_a^{-1}+e_1 \text{,}
\]
which is $0$.

The claim implies that the relators of the form $r=bc^{-1}db^{-1}cd^{-1}$ do not contribute in a significant way to the presentation of the $\Lambda$-module $M \otimes _{\Lambda_\mu}\Lambda$. For a relator $g_{a_1} g_{a_2} g_{a_1}^{-1} g_{a_3}^{-1}$ corresponding to a crossing, instead, the nonzero free derivatives are $1-g_{a_1} g_{a_2} g_{a_1}^{-1}$ with respect to $g_{a_1}$, $g_{a_1}$ with respect to $g_{a_2}$, and $-g_{a_1} g_{a_2} g_{a_1}^{-1} g_{a_3}^{-1}$ with respect to $ g_{a_3}$. The images in $\Lambda$ are $1-t$ with respect to $g_{a_1}$, $t$ with respect to $g_{a_2}$, and $-1$ with respect to $ g_{a_3}$. These coefficients provide the relation 
\[
0 = (1-t)((1 \otimes (g_{a_1}-1)) \otimes 1) +t((1 \otimes (g_{a_2}-1)) \otimes 1) - ((1 \otimes (g_{a_3}-1) ) \otimes 1)
\]
in $M \otimes _{\Lambda_\mu}\Lambda$, and this relation matches precisely with the crossing  relation for $\Mr(L)$ that appears in the definition given at the beginning of Sec.\ \ref{defs}.
\end{proof}

\begin{corollary}
\label{medonto}
If $L$ is a link with a diagram $D$, then there is a surjective quandle map $f_D:\QMG(L) \to Q^{\textup{red}}_A(L)$, under which $g_a \mapsto \gamma_D(a) \otimes 1 \thickspace \allowbreak \forall a \in A(D)$.  
\end{corollary}
\begin{proof}
The isomorphism $\widehat \gamma_D:M \otimes _{\Lambda_\mu}\Lambda \to \Mr(L)$ can be composed with the function $\MG(L) \to M \otimes_{\Lambda_\mu}\Lambda$ given by $g \mapsto (1 \otimes (g-1)) \otimes 1 \thickspace \allowbreak \forall g \in \MG(L)$. The function $f_D$ is obtained by restricting this composition to $\QMG(L)$. The composition certainly has $g_a \mapsto \gamma_D(a) \otimes 1 \thickspace \allowbreak \forall a \in A(D)$. The $\gamma_D(a) \otimes 1$ elements generate the quandle $\Qr(L)$, so in order to show that $f_D$ is a surjective quandle map, it is enough to show that $f_D$ is a quandle map.

If $x,y \in \QMG(L)$, then
\[
f_D(y \triangleright x) = f_D(xyx^{-1}) = \widehat \gamma_D((1 \otimes (xyx^{-1}-1))\otimes 1)
\]
\[
=\widehat \gamma_D((1 \otimes (xyx^{-1}-xy) )\otimes 1)+\widehat \gamma_D((1 \otimes (xy-x))\otimes 1)+\widehat \gamma_D((1 \otimes (x-1)) \otimes 1)
\]
\[
=\widehat \gamma_D((\alpha(xyx^{-1})\otimes (1-x) )\otimes 1)+\widehat \gamma_D((\alpha(x) \otimes (y-1))\otimes 1)+\widehat \gamma_D((1 \otimes (x-1)) \otimes 1)
\]
\[
=\widehat \gamma_D((1 \otimes (1-x) ) \otimes \tau \alpha(xyx^{-1}))+\widehat \gamma_D((1 \otimes (y-1)) \otimes \tau \alpha (x))+\widehat \gamma_D((1 \otimes (x-1)) \otimes 1).
\]

Every element of $\QMG(L)$ is a conjugate of some $g_a$ element, so its image under $\tau \alpha$ is $t$. Therefore $\tau \alpha (xyx^{-1}) = t^2t^{-1} = t = \tau \alpha(x)$, so
\[
f_D(y \triangleright x)=\widehat \gamma_D((1 \otimes (1-x) ) \otimes t)+\widehat \gamma_D((1 \otimes (y-1)) \otimes t)+\widehat \gamma_D((1 \otimes (x-1)) \otimes 1)
\]
\[
= t \cdot \widehat \gamma_D((1 \otimes (1-x) ) \otimes 1)+t \cdot \widehat \gamma_D((1 \otimes (y-1)) \otimes 1)+\widehat \gamma_D((1 \otimes (x-1)) \otimes 1).
\]
\[
= (1-t) \cdot \widehat \gamma_D((1 \otimes (x-1) ) \otimes 1)+t \cdot \widehat \gamma_D((1 \otimes (y-1)) \otimes 1).
\]
According to Proposition \ref{standardsub}, it follows that in $\Qr(L)$,
\[
f_D(y \triangleright x)=\widehat \gamma_D((1 \otimes (y-1) ) \otimes 1) \triangleright \widehat \gamma_D((1 \otimes (x-1) ) \otimes 1)=f_D(y) \triangleright f_D(x).
\]
 \end{proof}
\begin{corollary}
\label{disonto} Let $D$ be a diagram of $L$ with alternating writhes. Then the map $f_D$ induces an isomorphism $\Dis(f_D):\Dis(\QMG(L)) \allowbreak \to \Dis(\Qr(L))$ of abelian groups, given by the formula
\[
\Dis(f_D)(\prod \beta_{y_i}^{m_i}) = \prod \beta_{f_D(y_i)}^{m_i}.
\]
\end{corollary}
\begin{proof}
Corollary \ref{medonto} and property 6 of Proposition \ref{disprop} tell us that $f_D$ induces an epimorphism $\Dis(f_D)$, with the given formula. In order to show that $\Dis(f_D)$ is injective, we assemble an inverse function from pieces that have been discussed before. Let $a^* \in A(D)$ be fixed.  

According to Corollary \ref{subq}, $\Dis(Q_A^{\textup{red}}(L))$ is isomorphic to the $\Lambda$-submodule of $(1-t)M_A^{\textup{red}}(L)$ generated by $\{(1-t)(q-q') \mid q,q' \in Q_A^{\textup{red}}(L)\}$. This submodule is $(1-t) \ker \fr$. The isomorphism $g:\Dis(\Qr(L)) \to (1-t) \ker \fr$ given in the proof of Corollary \ref{subq} maps a displacement $\prod \beta_{q_i}^{m_i}$ to the module element $(1-t)(\sum m_i q_i)$. Therefore, if $a \in A(D)$ then this isomorphism $g$ has 
\[
g(\beta_{\gamma_D(a) \otimes 1} \beta_{\gamma_D(a^*) \otimes 1}^{-1})=(1-t)((\gamma_D(a) \otimes 1) - (\gamma_D(a^*) \otimes 1)) 
\]
\[
= (1-t)((\gamma_D(a) - \gamma_D(a^*)) \otimes 1).
\]

Let $e_D:\ker \fr \to \MG^0(L)$ be the $\Lambda$-linear map of Corollary \ref{dismap}, and let $\ann(1-t)= \{m \in \ker \fr \mid ((1-t)m=0 \}$. Proposition \ref{fed} tells us that $\ann(1-t)$ is the submodule of $\ker \fr$ generated by the longitudes $\chi_1, \dots \chi_\mu$.

An argument just like the proof of Proposition \ref{betamap} provides a homomorphism $\beta:\MG(L) \to \Aut(\QMG(L))$, with $\beta(g_a) = \beta_{g_a} \thickspace \allowbreak \forall a \in A(D)$. This homomorphism maps each product $\prod g_{a_i}^{m_i}$ to $\prod \beta_{g_{a_i}}^{m_i}$, so it maps $\MG^0(L)$ onto $\Dis(\QMG(L))$. According to Corollary \ref{dismod} and Theorem \ref{disqgen}, the scalar multiplication operations in the $\Lambda$-modules $\MG^0(L)$ and $\Dis(\QMG(L))$ can be defined using conjugation by $g_{a^*}$ and $\beta_{g_{a^*}}$, respectively. Clearly then the homomorphism $\beta$ is $\Lambda$-linear. 

An argument just like the proof of the third property of Proposition \ref{trivdis} shows that for each $i \in \{1. \dots, \mu\}$, the restriction of $ \beta e_D(\chi_i)$ to the $K_i$ orbit of $\QMG(L)$ is the identity map. Proposition 40 tells us that $\QMG(L)$ is semiregular, so it follows that each $\beta e_D(\chi_i)$ is the identity map of $\QMG(L)$. That is, $\chi_1, \dots, \chi_\mu \in \ker( \beta e_D)$. Then Proposition \ref{fed} implies that $\ann(1-t) \subseteq \ker ( \beta e_D)$, so $\beta e_D: \ker \fr \to \Dis(\QMG(L))$ induces a $\Lambda$-linear map $(\ker \fr) / \ann(1-t) \to \Dis(\QMG(L))$, given by $x+\ann(1-t) \mapsto \beta e_D(x)$. 

There is a $\Lambda$-linear epimorphism $\ker \fr \to (1-t) \cdot \ker \fr$ defined using scalar multiplication by $1-t$, and the kernel of this epimorphism is $\ann(1-t)$. Hence $(1-t) \cdot \ker \fr$ is isomorphic to $(\ker \fr) / \ann(1-t)$, and an isomorphism is given by $(1-t) \cdot x \mapsto x + \ann (1-t)$. Composing this isomorphism with the map induced by $\beta e_D$ that was mentioned at the end of the preceding paragraph, we obtain a map
\[
\widehat {\beta e_D}:(1-t) \cdot \ker \fr \to \Dis(\QMG(L)) \text{,}
\]
given by $\widehat {\beta e_D}((1-t)x) = \beta e_D(x)$
$\thickspace \allowbreak \forall x \in \ker \fr$. 

Recall the isomorphism $g$ mentioned in the second paragraph of the proof. The composition 
\[
\widehat {\beta e_D} \circ g: \Dis(\Qr(L)) \to \Dis(\QMG(L))
\]
is a $\Lambda$-linear map, and for every $a \in A(D)$ it has
\[
(\widehat {\beta e_D} \circ g)(\Dis(f_D)(\beta_{g_a}\beta_{g_{a^*}}^{-1}))=(\widehat {\beta e_D} \circ g)(\beta_{f_D(g_a)}\beta_{f_D(g_{a^*})}^{-1})
\]
\[
= \widehat {\beta e_D}(g(\beta_{\gamma_D(a) \otimes 1} \beta_{\gamma_D(a^*) \otimes 1}^{-1})) = \widehat {\beta e_D}((1-t)((\gamma_D(a) - \gamma_D(a^*)) \otimes 1))
\]
\[
= \beta e_D(((\gamma_D(a) - \gamma_D(a^*)) \otimes 1)) = \beta(h_a) = \beta(g_a g_{a^*}^{-1}) = \beta_{g_a}\beta_{g_{a^*}}^{-1}.
\]

The elementary displacements $\beta_{g_a}\beta_{g_{a^*}}^{-1}$ generate $\Dis(\QMG(L))$, so the composition $\widehat {\beta e_D} \circ g \circ \Dis(f_D)$ is the identity map of $\Dis(\QMG(L))$. As $\Dis(f_D)$ is surjective, it follows that $\Dis(f_D)$ is an isomorphism. \end{proof}

\begin{theorem}
\label{main3}
If $D$ is a diagram of $L$ with alternating writhes, then the quandle map $f_D:\QMG(L) \to \Qr(L)$ is an isomorphism.
\end{theorem}
\begin{proof}
Suppose $x \neq y \in \QMG(L)$, and $f_D(x)=f_D(y)$. 

According to Definition \ref{qmg}, there are arcs $a,b \in A(D)$ such that $x$ is a conjugate of $g_a$ and $y$ is a conjugate of $g_b$. Then the orbit of $\Qr(L)$ containing $f_D(x)=f_D(y)$ contains both $\gamma_D(a) \otimes 1$ and $\gamma_D(b) \otimes 1$, so $\kappa_D(a)=\kappa_D(b)$. Considering the crossing relations that appear in the description of $\MG(L)$ in Definition \ref{mg}, we deduce that $g_a$ and $g_b$ are conjugates of each other. Therefore $x$ and $y$ belong to the same orbit of $\QMG(L)$.

According to Proposition \ref{orb}, $y=d(x)$ for some $d \in \Dis(\QMG(L))$. As $x \neq y$, $d$ is not the identity map of $\QMG(L)$. According to Corollary \ref{disonto}, it follows that $\Dis(f_D)(d)$ is not the identity map of $Q_A^{\textup{red}}(L)$. As $Q_A^{\textup{red}}(L)$ is semiregular, this implies that $(\Dis(f_D)(d))(f_D(x)) \neq f_D(x)$. However, if $d=\prod \beta_{y_i}^{m_i}$ then using the fact that $f_D$ is a quandle map, we calculate
\[
(\Dis(f_D)(d))(f_D(x)) = \left( \prod \beta_{f_D(y_i)}^{m_i} \right) (f_D(x)) 
\]
\[
= f_D \left( \left( \prod \beta_{y_i}^{m_i}\right) (x) \right) = f_D (d(x))= f_D(y) = f_D(x) \text{,}
\]
a contradiction. We conclude that $f_D$ is injective. \end{proof}

\end{document}